\documentclass[11 pt]{amsart}
\usepackage{graphicx}

 \def\e{\hbox{\bf e}} 
 \def\x{\hbox{\bf x}} \def\y{\hbox{\bf y}}
 \def\z{\hbox{\bf z}}

 \def\sx{\hbox{\bf{\tiny x}}}
\def\sy{\hbox{\bf{\tiny y}}}

\def\nn{{\mathbb N}}  \def\rr{{\mathbb R}}

\def\I{\hbox{\bf I}}
\def\J{\hbox{\bf J}}
\def\R{\hbox{\bf R}}
\def\S{\hbox{\bf S}}
\def\sI{\hbox{\bf{\tiny I}}}
\def\sJ{\hbox{\bf{\tiny J}}}

\def\bs{\bigskip\noindent}
\def\ms{\medskip\noindent}

\def\AA{{\mathcal A}}
 \def\CC{{\mathcal C}} \def\DD{{\mathcal D}}
 \def\FF{{\mathcal F}} 
 \def\II{{\mathcal I}}
 \def\MM{{\mathcal M}} 
\def\OO{{\mathcal O}} \def\PP{{\mathcal P}} \def\QQ{{\mathcal Q}}
\def\RR{{\mathcal R}} \def\SS{{\mathcal S}}

\newtheorem{tmma}{\em Theorem}[section]
\newtheorem{lemma}{\em Lemma}[section]

\newtheorem{proposition}{\em Proposition}[section]
\newtheorem{corollary}{\em Corollary}[section]
\newtheorem{remark}{\em Remark}[section]

\def\endpf{\hfill{$\Box$}}

\begin{document}

%%%%%%%%%%%%%%%%%%%%%%%%%%%%%%%%%%%%%%%%

\baselineskip 16pt

\title[Complexity of regulatory networks ]{Dynamical complexity of 
discrete time regulatory networks}
\author{Ricardo Lima${}^1$ and Edgardo Ugalde${}^2$}
\maketitle

\bs 
\begin{center}
${}^1$Centre de Physique Th\'eorique \\
CNRS Luminy Case 907\\13288  Marseille, France

\ms

${}^2$ Instituto de F\'\i sica\\
Universidad Aut\'onoma de San Luis Potos\'\i \\
78000 San Luis Potos\'\i , M\'exico

\end{center}

\bs

\section{Introduction}

\

\noindent
The functioning and development of living organisms -- from bacteria to humans -- 
are controlled by genetic regulatory networks composed of interactions 
between DNA, RNA, proteins, and small molecules. The size and complexity 
of these networks make an intuitive understanding of their dynamics 
difficult to attain. In order to predict the behavior of regulatory systems 
in a systematic way, we need modeling and simulation tools with a solid 
foundation in mathematics and computer science.

\ms  In spite of involving a few number of genes, regulatory network
usually have complex interaction graphs. Many genes are connected to
each other and the incoming and outgoing interactions highly depend on
the gene. This is a major difference with spatially extended
systems where a usual assumption is that the incoming and outgoing
interaction do not depend on the node (translation
invariance). As in neural networks models, simplifying assumptions
are all--to--all interactions (globally coupled systems) or
random interactions. These assumptions seem not to fit exactly with
the interaction graphs constructed from biological data
analysis. Specific graph models and mathematical tools for the 
analysis of such complex networks of intermediate size are still
largely missing.

\ms Genes are segments of DNA involved in the production of proteins.
The process in which a protein is produced from the information encoded in the 
genes is called gene expression. During this process the genes are transcribed
into mRNA by enzymes called RNA polymerase, then the resulting mRNA 
functions as a template for the synthesis of a protein by another enzyme, 
the ribosome, in a process called translation. The level of
gene expression depends on the relative activity of the protein synthesis and 
degradation. In order to adapt the gene expression level to the requirements of
the living cell, the evolution has selected complicated mechanisms regulating the
production and degradation of proteins. A simple example of this, is the control 
of transcription by a repressor protein binding to a regulatory site on the DNA, 
preventing in this way the transcription of the genes. Hence, the expression 
level of a given gene regulates the expression level of another one, giving
rise to genetic regulatory system structured by networks of interactions between
genes, proteins, and molecules.  

\ms  Let us illustrate the type of interactions used 
to define our model by means of a simple example: the regulation 
of the expression of the sigma factor $\sigma^{S}$ 
in {\it Escherichia coli}, encoded by the gene {\it rpoS}. 
The factor $\sigma^{S}$ takes its name from the fact that it plays an
important role in the adaptation to a particular stress frequently
encountered by bacteria: the depletion of nutriments in the
environment, which leads to a considerable slowing down of cell
growth called the stationary growth phase. 
The protein CRP, a typical repressor--activator,
specifically binds the DNA at two sites close to the major promoter of
{\it rpoS}~\cite{HdJ2164}. One of these sites overlaps with the
promoter, which implies that CRP and RNA polymerase cannot
simultaneously bind to the DNA. As a
consequence, CRP represses the transcription of {\it rpoS}.
The example illustrates one type of regulation that is 
quite widespread in bacteria: a protein binding the DNA, the regulator, 
prevents or favors the binding of RNA polymerase to the promoter.

\ms  Genetic regulatory networks are usually 
modeled by systems of coupled differential equations~\cite{dejong}, 
and more particularly by systems of piecewise affine differential 
equations like in~\cite{edwards}.
Finite state models, better known as logical networks, are also used  (see~\cite{thieffrythomas}, or~\cite{thomaskaufman}). In this
paper we consider a class of models of regulatory networks which may
be situated in the middle of the spectrum; they present 
both discrete and continuous aspects. Our models 
consist of a network of units, whose states are quantified by 
a continuous real variable. The state of
each unit in the network evolves according to a contractive
transformation chosen from a finite collection of possible
transformations. Which particular transformation is chosen at a given 
time step depends on the state of the neighboring units.
These kind of models are directly inspired by the systems 
of piecewise linear differential equations mentioned before, 
though they do not correspond to a
time discretization of those differential equations, but rather to
a natural discrete time version of them. Such a discrete time model
is pertinent when delays in the communications between genes 
justifies discrete time evolution. Besides this, the great advantage 
of these kind of models is that they are suitable for a detailed 
study from the point of view of the theory of dynamical systems. 
In this framework we may formulate and answer general questions 
about the qualitative behavior of such a system, and eventually relate 
those results to the experiment.

\ms  A detailed description of the asymptotic dynamics of our models 
would require the use of symbolic dynamics, as it has been done 
in~\cite{gambaudo} for systems related to ours.
As a first approximation to the complete description of the dynamics  
we will focus on a global characteristics: the dynamical complexity. 
This is a well--studied notion in the framework of the theory of dynamical systems 
(see~\cite{ferenczi}, and~\cite{blanchardetal}), and it is related
to the proliferation of distinguishable temporal behaviors. The main motivation of
this work is to find explicit relations between the topological structure of the
regulatory network, and the growth rate of the dynamical complexity. More generally
we are interested in the possible constraints imposed 
by the topology of the underlying network and the contractive nature of the local 
dynamics, on the asymptotic dynamical behavior of such kind of systems. 
Hopefully, results on this direction will contribute to the
understanding of the observed behaviors in the biological systems inspiring our models.
The firs steps in this direction has been fulfilled by 
Theorems~\ref{polynomial}, \ref{skew} and~\ref{bound-degree} below.

\ms  In order to illustrate our results we show some biological motivated 
examples. Indeed, in subsection~\ref{ss-self-inhibitor} we discuss the self-inhibitor, 
whose biological prototype is presented in~\cite{self-inhibitor}. Then, in 
section~\label{s-examples} we apply our results to particular families of regulatory 
networks. First we consider the family of circuits which models concrete 
biological examples. For instance, the genetic toggle switch studied in~\cite{Gardner2000} 
corresponds to a circuit on two vertices, while the 
repressilator in {\it Escherichia coli}~\cite{Elowitz2000}, is modeled by a circuit on three 
vertices. Then, in a biological inspired four vertices network studied in~\cite{Vilela2004}, 
using the results of section~\ref{underlying-network}, we show that the complexity of 
this particular network may correspond to that of a network on two vertices.

\ms\noindent The paper is organised as follows. In the next section we
introduce the model, we give the basic definitions which we will
use in the sequel, and we present some basic facts concerning the
relation between complexity and asymptotic dynamics. In 
Section 3 we state and prove the main result of this paper, which concerns
the complexity of a discrete time regulatory network satisfying a 
(natural but) technical condition which we call coordinatewise
injectivity. In Section 4 we explore the relationship between the 
complexity and the structure of the network. Section 5 is devoted to 
some examples, and in the last section we present some concluding remarks.

\ms
\paragraph{\bf Acknowledgments} During the development of this
work, E. U. was invited researcher at the {\it Centre de Physique
Th\'eorique de Marseille}, and benefited from the financial
support of the {\it Minist\`ere d\'el\'eg\'ue \`a la recherche} through the 
{\it Universit\'e de Toulon et du Var}. Both
authors thank Bastien Fernandez for his remarks, comments, and
suggestions 
throughout %
this work.

\bs

\newpage
\section{Preliminaries}
\

\ms To each unit in interaction we associate a label in the finite
set $\{1,2,\ldots,d\}$. The interaction between units is codified
by an {\it interaction matrix} $K\in\MM_{d\times d}([0,1])$. The
matrix entry $K_{i,j}$ specifies the strength of the action of the
unit $i$ over the unit $j$. The matrix $K$ is normalized so 
that, %
for each $j$, $\sum_{i=1}^d K_{i,j}=1$. The interaction has two
possible modes: activation or inhibition. The interaction mode
between units is codified by an {\it activation matrix}
$s\in\MM_{d\times d}(\{-1,0,1\})$, which is compatible with the
interaction matrix, i.~e., $s_{i,j}=0$ if and only if $K_{i,j}=0$.
The strength of the interaction is weighted by a {\it contraction
rate} $a\in [0,1]$; it is directly related to 
the %
speed of degradation
of the state of a noninteracting unit. At each time step
$t\in\nn$, the state of a unit $j\in \{1,2,\ldots,d\}$ is
quantified by a real value $\x^t_j\in [0,1]$, and it changes in
time according to the transformation
\[
\x_j^{t+1}=a\x_j^t+(1-a)\sum_{i=1}^dK_{i,j}H(s_{i,j}(\x_i^t-T_{i,j})),
\]
where $H:\rr\to\{0,1\}$ is the Heaviside function,
\[
x\mapsto H(x):=\left\{ \begin{array}{ll}
                       0 & \text{ if } x\leq 0\\
                       1 & \text{ otherwise.}\end{array}
                             \right.
\]
All discontinuities of the transformation are contained in the
{\it threshold matrix} $T\in\MM_{d\times d}([0,1])$, which is
compatible with the interaction matrix.

\ms We let %
 the mapping $F:[0,1]^d\rightarrow[0,1]^d$ defined by 
\begin{equation}\label{map}
F(\x)_j:=
a\x_j^t+(1-a)\sum_{i=1}^dK_{i,j}H(s_{i,j}(\x_i^t-T_{i,j}))
\end{equation}
induce %
a discrete dynamical system on the $d$--dimensional unit
cube, a {\it discrete time regulatory network}. Due to
discontinuities, this dynamical system may exhibit a complex
behavior characterized by the existence of several attractors
whose basins of attraction are intermingled in a complicated way.
Often these attractors consist of a collection of disjoint
periodic orbits, but there 
also exist %
attractors with infinite cardinality.

\ms
For $a = 0$, the discrete time regulatory network is equivalent to a 
{\it logical network}~\cite{glass, thomas}.
For $a > 0$, the model is no longer equivalent to a boolean network, 
and can be viewed as a discrete time analogue of the a system of delay 
differential equations (see~\cite{bastien} for details).
The contraction rate $a$ plays the role of (the inverse of) a 
delay parameter: the smaller $a$, the stronger the delay is.
\ms

\subsection{The self--inhibitor}~\label{ss-self-inhibitor}
In~\cite{self-inhibitor} it is studied an autoregulatory system which
was experimentally implemented to investigate the role of negative
feedback loops in the gain of stability in genetic networks. This 
system can be modeled by the one--dimensional regulatory network
 		$x\mapsto F(x):=a x + (1-a)H(T-x)$, which
is semiconjugated to a rotation~\cite{coutinho}. As we move the
parameters $a$ and $T$,  the associated rotation number ranges
over an interval, so that the attractor is either a finite set, in
the case of a rational rotation number, or a Cantor
set in the irrational case. Both regimes of the self--inhibited
unit, periodic and quasiperiodic, share a common feature: the
number of distinguishable orbits grows at most linearly in time.
To be more precise, the orbits $\{x^0, x^1=F(x^0),\ldots,
x^t=F(x^{t-1}),\ldots\}$ and  $\{y^0, y^1=F(x^0),\ldots,
x^t=F(y^{t-1}),\ldots\}$ are distinguishable at time $t=\tau$ if
there exists $0\leq t\leq \tau$ such that $x^t$ and $y^t$ lie at
opposite sides of the discontinuity, i.~e., $H(T-x^t)=1-H(T-y^t)$.
Note that due to the piecewise contractive nature of the dynamics,
two indistinguishable initial conditions get exponentially close
as time goes to infinity.

\ms The asymptotic dynamics of the self--inhibited unit takes
place inside the interval $I:=[aT, aT+(1-a)]$, where every initial
condition get trapped after a finite number of iterations. Inside the attractor
there are two classes of initial 
conditions, % 
corresponding to the
intervals $I_0=[T,aT+(1-a)]$ and $I_1:=[aT,T)$. Since the map
$x\mapsto F(x):=ax+(1-a)H(T-x)$ is injective inside $I$, and it is
an affine contraction when restricted to $I_0$ or to $I_1$, then
at most one of the intervals $F(I_0)$, $F(I_1)$, contains the
discontinuity. The other one would be completely contained in one
of the original intervals $I_0$ or $I_1$. Thus, the interval whose
image lies inside one of the original intervals contains initial
conditions whose orbits are indistinguishable at time $t=1$,
whereas the other interval contains initial conditions which may
generate up to two distinguishable orbits. Therefore, for $t=1$
there are at most 3 distinguishable orbits represented by the same
number of disjoint intervals. Those intervals are precisely the
domains of continuity of the mapping $x\mapsto F^2(x)$. We may
pursue this reasoning and show by induction that for each $t\in
\nn$, the maximal cardinality of a set of initial conditions in
$I$ which generate mutually distinguishable orbits at time $t$, is
at most $t+2$.

\ms For the values of $a$ and $T$ for which the dynamics is
semiconjugated to an irrational rotation, the number of
distinguishable orbits may grow indefinitely. Indeed, for those
values, for each $n\in\nn$ there exists one interval $J\subset I$
where $F^n$ is a contraction, and such that $T\in F^t(J)$. This is
why for these %
particular choices of $a$ and $T$ % 
there are $t+1$
domains of continuity of $F^t$, for each $t\in\nn$. This is
precisely the case when the asymptotic dynamics is quasiperiodic,
semiconjugated to an irrational rotation (see~\cite{coutinho} for
details).

\ms In the sequel of this paper we will extend these arguments to
the general case, and in this way we will obtain upper bounds for
the number of distinguishable orbits.

\ms Let us remark that the system studied in~\cite{gambaudo} is a generalization
of the self--inhibitor we present here.

\ms

\subsection{The complexity} In the general $d$--dimensional case,
given a threshold matrix $T\in \MM_{d\times d}([0,1])$, for each
coordinate $1\leq i\leq d$ define the partition $\PP_i$ of $[0,1]$
so that the map
\[
\x_i\mapsto ( H(s_{i,1}(\x_i-T_{i,1})),
H(s_{i,2}(\x_i-T_{i,2})),\cdots, H(s_{i,d}(\x_i-T_{i,d})) )
\]
is constant restricted to each atom of $\PP_i$. We can see that
the atoms of $\PP_i$ are semiclosed intervals with endpoints in
the set $\{T_{i,j}: \ j=1\ldots,d\}\cup\{0,1\}$. The map
$F:[0,1]^d\to [0,1]^d$ defined by Eq.~(\ref{map}) %
is an
affine contraction restricted to the atoms of
the %
Cartesian %
product
$\PP:=\prod_{i=1}^d\PP_i$, which we call {\it base partition}.
Indeed, if $\x,\y\in [0,1]^d$ belong to the same atom of $\PP$,
then for each $1\leq i,j\leq d$ we  have
$H(s_{i,j}(\x_i-T_{i,j}))=H(s_{i,j}(\y_i-T_{i,j}))$, which implies
$F(\x)-F(\y)=a(\x-\y)$.

\ms We define recursively the dynamical partitions
\begin{equation}\label{partitions}
\PP^{t+1}:=\{F^{-1}(\I)\cap \J:\ \I\in \PP^t \text{ and } \J\in
\PP\},
\end{equation}
with $\PP^1:=\PP$. It is easy to verify that the $\tau$--th
iterate of the map $F$, $F^{\tau}:[0,1]^d\to[0,1]^d$, is an affine
contraction in each atom of $\PP^{\tau}$. It means that two
initial conditions in the same atom of $\PP^{\tau}$ will 
get closer during the first $\tau$ time steps.

\ms The following discussion holds in the slightly more general
case when the mapping $F:[0,1]^d\to[0,1]^d$ is a non necessarily affine 
contraction when
restricted to each atom of the base partition $\PP$. We call this
{\it a piecewise contraction}.

\ms {\it The complexity} of the system $([0,1]^d,F)$ is the
integer valued function
\[
C:\nn\to\nn \ \text{ such that}\ C(t):=\#\PP^t.
\]
The number of atoms in the dynamical partition $\PP^{t}$ equals
the number of distinguishable orbits up to time $t$. It is the aim
of this paper to give upper bounds on the growth of this quantity. 
For systems defined by a continuous transformation
acting on a compact metric space, if the complexity is eventually
constant then all orbits are eventually periodic. For the same
kind of system, to a fast growth of the complexity one can
associate some notion of chaos (see \cite{blanchardetal} and
\cite{ferenczi}). The computation of the complexity is part of
program whose aim is the complete characterization of the
asymptotic dynamics of the discrete time regulatory network. 
Information about the growth of the complexity,
together with a description of the recurrence properties of
individual orbits, will allow us to give a detailed picture of the
asymptotic behavior of the system.

\ms \subsection{The asymptotic behavior} {\it The limit set}
$\Lambda$ of the system $([0,1]^d,F)$ is the collection of all
accumulation points of its orbits, i.~e., $\x\in\Lambda$ whenever
there exists an initial condition $\y\in [0,1]^d$, and a sequence
of times $t_1 < t_2 < \cdots <t_n <\cdots $, such that
$\lim_{n\to\infty}F^{(t_n)}(\y)=\x$.

\ms Let us denote $$\Delta:=\left\{\x\in [0,1]^d:\ \exists i,
k\in\{1,2,\ldots, d\}\ \text{such that}\  \x_k=T_{i,k}\right\},$$
the discontinuity set. When $\Lambda \cap \Delta = \emptyset$,
then $\Lambda:=\bigcap_{t=0}^\infty F^t([0,1])^d$, and it 
contains all the %
orbits with infinite past.

\ms If the limit set has infinite cardinality, then the complexity
is a strictly increasing function. Indeed, if on the contrary
$C(\tau+1)=C(\tau)$ for some $\tau\in\nn$, then necessarily to
each atom $\I\in\PP^{\tau}$ 
there %
corresponds a unique atom
$\bar{\I}$ such that $F(\I)\subset \bar{\I}$. The transformation
$\I\mapsto \bar{F}(\I):=\bar{\I}$, acting over $\PP^{\tau}$, would
generate in this case up to $C(\tau)$ ultimately periodic orbits.
Thus, to each atom $\I\in \PP^{\tau}$ we could associate a
transient $t\in \nn$ and a period $p\in \nn$ such that
$\bar{F}^{t+kp}(\I)=\bar{F}^t(\I)$ for each $k\in \nn$. Since $F$
is a contraction when restricted to each atom of $\PP^{\tau}$, all
the orbits with initial condition $\x\in\I$ would have one and the
same accumulation point in the closure of each one of the atoms
$\bar{F}^t(\I), \bar{F}^{t+1}(\I),\ldots, \bar{F}^{t+p-1}(\I)$.
Therefore, $C(\tau+1)=C(\tau)$ implies $\#\Lambda < \infty$.
Of course, if $\Lambda$ is a finite set, then the dynamical
complexity has to be a bounded function.

\ms One would expect bounded dynamical complexity to be a generic
situation in the case of piecewise contractions. 
The following suggests this is true.  %%%%%

\ms
\begin{proposition}[Multiperiodicity]
If ${\rm dist}(\Lambda,\Delta)>0$ then $\# \Lambda < \infty$.
\end{proposition}

\ms
\begin{proof}
Let $\epsilon={\rm dist}(\Lambda,\Delta)/2> 0$. Compactness allows us to
choose a finite cover $\CC=\{B(\x^{(i)},\epsilon_i):\ \x^{(i)}\in
\Lambda,\ i=1,2,\ldots N\}$ for $\bar\Lambda$, composed by open
balls of radius not greater than $\epsilon$ and centered at points
in $\Lambda$. Since $\epsilon_i < {\rm dist}(\Lambda,\Delta)$ for
each $i=1,2,\ldots, N$, then $F^t(B(\x^{(i)},\epsilon_i))$ is a
ball of radius $a^t\epsilon_i$ centered at $F^t(\x^{(i)})$, for
each $t\in\nn$. In this way, for each $t\in \nn$, the collection
$\CC_t:=\{B(F^t(\x^{(i)}), a^t\epsilon_i):\ i=1,2,\ldots,N\}$ is
an open cover for $F^t(\bar\Lambda)$, but since $\Lambda$ is
$F$--invariant, then $\CC_t$ is an open cover for $\Lambda$.

\ms Let $t$ be such that $a^t \epsilon$ is less than or equal to
the Lebesgue number of $\CC$. Then, for each $i=1,2,\ldots, N$
there exists $j\in \{1,2,\ldots,N\}$ such that
$$B(F^t(\x^{(i)}),a^t \epsilon_i)\subset B(\x^{(j)},\epsilon_j).$$
For each $i=1,2,\ldots,N$ choose one such $j$, and define a map
$i\mapsto \tau(i):=j$ from $\{1,2,\ldots,N\}$ to itself. The
successive iterations of this map, $i\mapsto \tau(i)\mapsto
\cdots$, generate ultimately periodic sequences, but the
$F$--invariance of $\Lambda$ implies that each $i\in
\{1,2,\ldots,N\}$ has to be periodic under $\tau$. This means that
for each $i\in\{1,2,\ldots,N\}$ there exists $p=p(i)\in \nn$ such
that
\[
F^{m\times tp}(B(\x^{(i)},\epsilon_i))\subset
B(\x^{(i)},\epsilon_i)\ \forall \ m\in \nn.
\]
Contractivity implies that $B(\x^{(i)},\epsilon_i)\cap
\Lambda=\{\x^{(i)}\}$ for each $i\in \{1,2,\ldots,N\}$.
Furthermore, $\x^{(i)}$ is periodic of period (not necessarily
minimal) $t\times p(i)$. Periods are bounded by $N\times t$, so
they depend on ${\rm dist}(\Lambda,\Delta)$ through the Lebesgue
number of $\CC$.
\end{proof}

\ms We would like to prove that, in the case of
networks of piecewise affine contractions, ${\rm
dist}(\Lambda,\Delta)>0$ is generic.  This is true for in the 1--dimensional case~\cite{coutinho}, but the question remains open for $d\geq 2$.

\ms When on the contrary ${\rm dist}(\Lambda,\Delta)=0$, the
attractor may fail to be invariant under the dynamics. In this
case, to each point $\x\in \Lambda$ we may associate a ``ghost
orbits'', which would be a real orbit by a convenient
redefinition of the transformation at the discontinuities.

\bs Given an orbit $O(\x)$ and $\tau\in\nn$, let
$\PP^{\tau}_{\sx}:=\{\I\in\PP^{\tau}:\ \exists\ t\in\nn \text{
such that } \x^{(t)}\in\I\}$. This is the restriction of the
dynamical partition to the orbit. We will say that the orbits
$O(\x)$ and $O(\y)$ are distinct if $\PP^{\tau}_{\sx}\cap
\PP^{\tau}_{\sy}=\emptyset$ for some $\tau\in \nn$. Note that
distinct orbits are necessarily distinguishable at some time.
i.~e., if $O(\x)$ and $O(\y)$ are distinct then there exists
$\tau\in\nn$ such that $\x^{\tau}$ and $\y^{\tau}$ belong to
different atoms of the base partition $\PP$.

\ms We can easily relate the growth of the dynamical complexity to
the number of infinite, distinct orbits. Indeed, for $O(\x)$
infinite, we have $\#\PP^t_{\sx}\geq t+1$ for all $t\in\nn$.
Otherwise there would exist $t\in\nn$ such that $\PP^t=\PP^{t+1}$ and
the orbit would be ultimately periodic. Hence, if $\Lambda$
contains $K$ distinct infinite orbits $O(\x^{1}), O(\x^{2}),
\ldots, O(\x^{K})$, 
then there exists a time $\tau\in \nn$ such that
$C(t)\geq \sum_{i=1}^K \#\PP^t_{\sx^i}\geq K(t+1)$ for each $t\geq
\tau$. Therefore, assuming the lowest possible complexity for each
infinite component of the attractor, a linear behavior of the
complexity would imply the finiteness of the number of such
components.

\bs The dynamical complexity is directly related to the
topological entropy. For a dynamical system $(X, F)$, for each
$\tau\in \nn$ and $\epsilon
>0$, a set $E\subset X$ is said to be 
$(\tau,\epsilon)$--{\it separated} if for $\x,\y\in E$ there exists 
$0\leq t\leq \tau$ such that $d(F^t(\x),F^t(\y))\geq \epsilon$. The
topological entropy of the system $(X, F)$ is defined by
\[h_{\rm top}:=\lim_{\epsilon\to 0}(\limsup_{\tau\to\infty}\log
N(\tau,\epsilon)/\tau),
\]
where $N(\tau,\epsilon)= \max\{\#E:\ E\subset X \text{ is }
(\tau,\epsilon)\text{--separated}\}$.

\ms The topological entropy of a piecewise contraction $(X,F)$
satisfies the inequality
\[h_{\rm top}\leq \limsup_{t\to\infty}\frac{1}{t}\log C(t).\]
Indeed, for $\tau\in\nn$, $\epsilon > 0$ and $\I\in \PP^{\tau}$
fixed, if $E\subset X$ is $(\tau,\epsilon)$--separated, then
$E\cap \I$ is necessarily $(1,\epsilon)$--separated. Therefore
$N(\tau,\epsilon)\leq C(\tau)\times N(1,\epsilon)$.

\section{Main Result}

\ms
\begin{tmma}[Polynomial upper bound]\label{polynomial}
Given the interaction and threshold matrices $K, T\in \MM_{d\times
d}([0,1])$, and the activation matrix $s\in\MM_{d\times
d}(\{-1,0,1\})$, there exist  $a_0\in [0,1]$ and $c\geq 1$, such
that for each contraction rate $a\in [0, a_0)$, the complexity of
corresponding discrete time regulatory network satisfies $C(t)
\leq 1+ c \left(1+c^d t^d\right)$.
\end{tmma}

\ms In the next section we will consider some specific contraints on the
structure of the network with which we are able to improve this upper
bound. We will also show some concrete examples where the constant 
involved in the upper bound is explicitly determined.

\ms 

\subsection{Prerequisites} 
For each $1\leq j\leq d$ let us define the system of affine
1--dimensional contractions
\begin{equation}\label{syslocal}
\FF_j:=\left\{x\mapsto
ax+(1-a)\left(\sum_{i=1}^d\e_iK_{i,j}\right) :\ \e\in\{0,1\}^d
\right\}.
\end{equation}
This system of contractions allows us to compute the coordinates
of $F$. Indeed, for each $\I\in \PP$ and $1\leq j\leq d$, there
exists an affine contraction $f\in \FF_i$ such that
$F(\x)_j=f(\x_j)$ whenever $\x\in \I$.

\ms The discrete time regulatory network $([0,1]^d,F)$ is
said to satisfy {\it coordinatewise injectivity} if for
each $1\leq j\leq d$
and $f,f'\in \FF_j$ we have
\[
f([0,1])\cap f'([0,1])=\emptyset\ \text{ whenever }\ f\neq f'.
\]
To prove the main result of this paper we first show that for
networks satisfying this property, $C(t)$ grows at most
polynomially. In particular, for a network satisfying
coordinatewise injectivity, the corresponding mapping $F$ is
injective. In this case the sets in $\QQ^{t+1}:=F^t(\PP^{t+1})$
are disjoint affine images of the atoms in $\PP^t$, and we can use
the equation $\#\PP^t=\#\QQ^t$ to compute the complexity. This
equality holds in the more general situation, when for each $t\in
\nn$, $F^t$ maps different atoms in $\PP$ to different sets, not
necessarily disjoint.

\ms In the coordinatewise injective case, the degree of the
polynomial bound for the complexity may be strictly smaller
that the cardinality of the network. In the next section we will
see how this degree can be related to the topological structure of
the network.

\bs \begin{lemma}~\label{injective} If $([0,1]^d,F)$ satisfies
coordinatewise injectivity, then %
$F$ is injective.
\end{lemma}

\begin{proof}
For $\x,\y \in [0,1]^d$ with $\x\neq \y$ we have two
possibilities: either they belong the same atom of
the base partition, or  they belong to different atoms.

\ms In the first case, since the piecewise constant 
part of $F$ 
\[
\x \mapsto (1-a)\left( \sum_{i=1}^d K_{i,j} H(
s_{i,j}(\x_i-T_{i,j}) )\right)_{j=1}^d 
\]
is constant on the atoms of $\PP$, 
it follows that $F(\x)-F(\y)=a(\x-\y)\neq 0$.

\ms In the second case there exists $1\leq i\leq d$, and different
affine contractions $f,f'\in \FF_i$, such that $F(\x)_i=f(\x_i)$ while 
$F(\y)_i=f'(\y_i)$. The coordinatewise
injectivity ensures that $f$ and $f'$ have disjoint images, which
implies that $F(\x)_i\neq F(\y)_j$, 
and hence $F(\x)-F(\y)\neq 0$ as required. 
\end{proof}

\ms
\begin{proposition}[Strong contraction implies coordinatewise 
injectivity]~\label{prop-injectivity} For each
$K\in\MM_{d\times d}([0,1])$ there exists $a_0 > 0$ such that any
discrete time regulatory network with interaction matrix $K$ and
contraction rate $a\in[0,a_0)$ satisfies coordinatewise
injectivity.
\end{proposition}

\ms

\begin{proof} Fix the threshold and interaction matrices
compatible with $K$. For each $j\in \{1,2,\ldots,d\}$, let $\FF_j$
be the system of affine contractions defined in (\ref{syslocal}).
For $j\in \{1,2,\ldots,d\}$ fixed, let
\[
\II_j:=\left\{\eta\in [0,1]:\ \eta:=\sum_{i=1}^d\e_iK_{i,j} \text{
for some } \e\in \{0,1\}^d\right\}.
\]
Let $\delta_j:=\min\{|\eta-\eta'|:\ \eta,\eta'\in \II_j\text{ and
} \eta\neq \eta'\}$. Finally, let $\delta:=\min_{1\leq j\leq
d}\delta_j$, which is strictly positive. With this define
$a_0:=\delta/(1+\delta)$, so that $a<(1-a)\delta$ for each $0\leq
a<a_0$.

\ms Suppose that $a < a_0$, fix $i\in\{1,2,\ldots,d\}$, and take
$f,f'\in \FF_i$ such that $f\neq f'$. Hence there are
$\eta,\eta'\in \II_i$ with $\eta\neq \eta'$, such that
$f(x):=ax+(1-a)\eta$ and $f'(x)=ax+(1-a)\eta'$. Without loss % 
of
generality, let us suppose that $\eta > \eta'$, which implies that
$\eta \geq \eta'+\delta$. Then we have
\[
\max_{x\in [0,1]}f'(x)=a+(1-a)\eta' < (1-a)(\eta'+\delta) <
                                      (1-a)\eta=\min_{x\in [0,1]}f(x),
\]
which holds for $i\in \{1,2,\ldots,d\}$ arbitrary, and the
proposition follows.
\end{proof}

\ms For each $1\leq i\leq d$ and $t\in \nn$, let $\QQ_i^t$ be the
collection of intervals
\[
\QQ^t_i:=\{J\subset [0,1]:\ J=\Pi_i\J \text{ for some } \J
\in\QQ^t\},
\]
Let us recall that $\QQ^{t+1}$ is the image, at time $t$, of the
dynamical partition $\PP^{t+1}$. Here and below $\Pi_i$ will denote 
the projection on the $i$--th coordinate.

\bs \begin{lemma}~\label{mutually-disjoint} If $([0,1]^d,F)$
satisfies coordinatewise injectivity, then for each $1\leq
i\leq d$ and $t\in \nn$, the intervals in $\QQ_i^t$ are 
pairwise disjoint.
\end{lemma}

\ms \begin{proof} For each $1\leq i\leq d$ $\QQ_i^1=\PP_i$ is a
partition, and so the intervals in $\QQ_i^1$ are 
pairwise disjoint.
Suppose that for $t=\tau$ the intervals in $\QQ_i^{\tau}$ are
partwise disjoint. Let $\FF_i$ be the system of affine
contractions as defined in equation (\ref{syslocal}). Then 
\[\QQ_i^{\tau+1}\subset \FF_i(\QQ_i^{\tau})\vee\PP_i:=
\{ f(J)\cap I:\ f\in \FF_i, \ J\in \QQ_i^{\tau}, \text{ and } I\in
\PP_i\}.
\]
Since $F$ satisfies coordinatewise injectivity, and since the
intervals in $\QQ^{\tau}_i$ are mutually disjoint, it follows that 
$\FF_i(\QQ_i^{\tau})\vee\PP_i$ is a collection of 
pairwise disjoint intervals, implying that the intervals in
$\QQ_i^{\tau+1}$ 
are partwise disjoint as well.
\end{proof}

\ms Given a set of coordinates $U\subset \{1,2,\ldots,d\}$, and
$t\in\nn$ fixed, a {\it $(U, t)$--specification} is a
$\#U$--dimensional interval $\S\in \prod_{i\in U}\QQ^t_i$. The
{\it bouquet} corresponding to the specification $\S$ is the
collection of $d$--dimensional intervals
\begin{equation}~\label{bouquet}
C(\S):=\{\J\in \QQ^t:\ \Pi_i\J=\Pi_i\S \ \forall\ i\in V'\}.
\end{equation}
We will say that $\J\in \QQ^t$ satisfies specification $\S$ if
$\J\in C(\S)$, and we will call $N(\S):=\#C(\S)$ the {\it degeneracy} 
of $\S$.

\ms
\begin{corollary}~\label{determine}
If $([0,1]^d,F)$ satisfies coordinatewise injectivity, for each
$1\leq i, j\leq d$ and $t\in\nn$, there exists an interval
$J^t_{i,j}\in\QQ_i^t$ such that $\{\J\in\QQ^t:\
T_{i,j}\in\Pi_iF(\J)\}\subset C(J^t_{i,j})$.
\end{corollary}

\ms
\begin{proof} For each  $\J\in \QQ^t$ we have
$\Pi_i F(\J)\in \FF_i(\QQ^t_i):= \{f(J):\ f\in \FF_i \text{ and }
J\in \QQ^t_i\}$. The previous lemma ensures that $\FF_i(\QQ^t)$ 
consists of partwise disjoint intervals. Thus, for each $(i,j)\in\AA$,
there cannot be more that one interval $J\in \QQ^t$, and one
affine contraction $f\in \FF_i$, such that $T_{i,j}\in f(J)$.
\end{proof}

\ms We denote $J_{i,j}^t$ the unique interval in $\QQ_i^t$ such
that if $T_{i,j}\in \Pi_iF(\J)$ for $\J\in \QQ^t$, then
$\Pi_i\J=J_{i,j}^t$.

\ms We say that the $d$--dimensional interval $\J\in\QQ^t$ is a
{\it predecessor} of $\J'\in\QQ^{t+1}$ if $\J'=F(\J)\cap \I$ for
some $\I\in \PP$. Reciprocally we say that $\J'$ is a {\it
successor} of $\J$. In our case, since $F$ is injective, every
$\J'\in \QQ^{t+1}$ has a unique predecessor $\J\in\QQ^t$, which we
denote by $P(\J')$.

\begin{lemma}~\label{predecessor}
Suppose that $([0,1]^d,F)$ is coordinatewise injective. Fix
$U\subset \{1,2,\ldots,d\}$ and $t\in\nn$. For each
$(U,t+1)$--specification $\S$, there is a unique
$(U,t)$--specification $P(\S)$, such that $\J\in C(\S)$ implies
$P(\J)\in C(P(\S))$.
\end{lemma}

\begin{proof}
Fix $i\in U$ and $\J\in \QQ^{t+1}$ satisfying the specification
$\S$. Let $\J'=P(\J)$ and $\I\in \PP$ such that $\J=F\J'\cap \I$.
Then 
\begin{equation}\label{father}
\Pi_i\S=\Pi_i(F(\J')\cap \I)= f(\Pi_i\J')\cap\Pi_i\I 
\end{equation}
for some $f\in\FF_i$. Since $F$ is coordinatewise 
injective, then the
affine contraction $f\in\FF_i$ for which (\ref{father}) holds is
uniquely determined from $\J'\in\QQ^t$. On the other hand,
Lemma~\ref{mutually-disjoint} establishes that for each $1\leq i\leq d$ and
$t\in \nn$, the collection $\QQ^t_i$ 
consists of partwise disjoint
intervals. Hence, for $f\in\FF_i$ fixed, there is only one
interval $J'\in \QQ_i^t$ such that $\Pi_i\J'=J'$ satisfies
Equation (\ref{father}). Therefore $\Pi_i\J'$ depends only on $\Pi_i\S$,
and the lemma is proved.
\end{proof}

\ms
\begin{lemma}\label{recurre1}
Let $([0,1]^d, F)$ be coordinatewise injective, and let $U\subset
\{1,2,\ldots,d\}$. For each $t\in \nn$, and every
$(U,t+1)$--specification $\S$, we have
\[
N(\S)\leq N(P(\S)) +
c\times \sum_{i\not\in U}\sum_{j=1}^d N(P(\S)\times J_{i,j}^t),
\]
with $J_j^t$ as defined by above, and $c:=\#\PP-1$.
\end{lemma}

\begin{proof}
Lemma \ref{predecessor} implies that
\[
C(\S)=\{F(\J')\cap\I\in
C(\S):\ \I\in \PP \text{ and } \J'\in C(P(\S)) \}.
\]
If $\J'\in C(P(\S))$ has more than one successor
satisfying $\S$, then necessarily $T_{i,j}\in \Pi_iF(\J')$ for
some $i\not\in U$ and $1\leq j\leq d$. Now, each $\J'\in C(P(\S))$
has no more than $\#\PP$ successors. Therefore %
\[
N(\S)\leq N(P(\S))+c \times \#\{\J'\in C(P(\S)):\
T_{i,j}\in\Pi_iF(\J')\text{ for some }i\not\in U\text{ and}
1\leq j\leq d\}.
\]

\ms From Lemma \ref{determine}, if $\J'\in C(P(\S))$ has more than
one successor 
that satisfies % 
$\S$, then $\J'$ satisfies
$P(\S)\times J_{i,j}^t$ for some $i\not\in U$ and $1\leq j\leq d$.
Thus,
\[
\{\J'\in C(P(\S)): T_{i,j}\in \Pi_iF(\J')\text{ for some }i\not\in
U \text{ and } 1\leq j\leq d\} \subset \bigcup_{i\not\in
U}\bigcup_{j=1}^d C(P(\S)\times J_{i,j}^t),
\]
and from this we readily obtain the required
inequality.
\end{proof}

\bs 
\subsection{Proof of Theorem~\ref{polynomial}} Fix
$a\in[0,a_0)$, with $a_0>0$ depending on $K,T\in \MM_{d\times
d}([0,1])$ and $s\in\MM_{d\times d}(\{-1,0,1\})$ as in the proof
of Propositon~\ref{prop-injectivity}. This proposition ensures that the
mapping $F$ satisfies coordinatewise injectivity.
Lemma~\ref{injective} ensures that $F$ is injective, 
and so we
can replace $\PP^t$ by $\QQ^t$ in the computation of the
complexity. In order to keep control of the 
growth of $\#\QQ^t$, we use the recurrence $\QQ^{t+1}=F(\QQ^t)\vee \PP$.

\ms If $C(t+1) > C(t)$, then %
the additional contribution in $C(t+1)$ is
due to those $d$--dimensional intervals in $\QQ^t$ which posses
more than one successor.

\ms Since $\J\in\QQ^t$ cannot have more than $\#\PP$ successors in
$\QQ^{t+1}$, one for each possible intersection $F(\J)\cap \I$
with $\I\in\PP$, 
it follows that the complexity satisfies the recursion
\begin{equation}\label{recurrenceC}
C(t+1)\leq C(t)+c_1 \times \#\{\J\in \QQ^t:\ \J \text{ has more
than two successors}\},
\end{equation}
where $c:=\#\PP-1$.

\ms Now, for $\J\in\QQ^t$ to have more than one successor it is
required that $F(\J)$ intersect at least two different atoms in
$\PP$. But if this is the case, then necessarily $F(\J)$ contains
a discontinuity set for $F$. Hence, if $\J\in\QQ^t$ has % 
more than
one successor, then %
there exists $1\leq i,j\leq d$ such that
$T_{i,j}\in \Pi_iF(\J)$. Then, according to (\ref{recurrenceC}) we
have
\begin{equation}\label{intersects}
C(t+1)\leq C(t)+c\times \#\{\J\in\QQ^t: \ T_{i,j}\in\Pi_iF(\J)
\text{ for some } 1\leq i,j\leq d\}.
\end{equation}

\ms Using lemma \ref{determine} and % 
(\ref{intersects})
we conclude that
\[
C(t+1)\leq C(t)+c\times \#\{\J\in\QQ^{t}:\ \Pi_i\J=J^{t}_{i,j}
\text{ for some } 1\leq i,j\leq d\},
\]
and so we can rewrite~(\ref{intersects}) as
\begin{equation}\label{recurrenceC2}
C(t+1)\leq C(t)+c\times \sum_{i=1}^d\sum_{j=1}^d
\#\{\J\in\QQ^t:\ \Pi_i\J=J_{i,j}^t\}.
\end{equation}
Since the cardinalities
$\#\{\J\in\QQ^t:\ \Pi_j\J=J_{i,j}^t\}$ in (\ref{recurrenceC2}) are the
same as the degeneracies $N(J_{i,j}^t)$,
we can rewrite 
(\ref{recurrenceC2}) as
\begin{equation}\label{recurrenceC3}
C(t+1)\leq C(t)+c\times \sum_{i=1}^d\sum_{j=1}^d N(J_{i,j}^t),
\end{equation}
where $c:=\#\PP-1$.

\ms In what follows we will establish polynomial bounds for the
growth of the degeneracies $N(J)$, with $J\in \QQ_i^t$ arbitrary.
For this aim we use a recurrence relation between degeneracies of
successive specifications. According to Lemma~\ref{predecessor},
for $U\subset \{1,2,\ldots,d\}$ and $t\in\nn$ fixed, to each
$(U,t+1)$--specification $\S$, there %
corresponds a unique
$(U,t)$--specification $P(\S)$, such that if $\J\in\QQ^{t+1}$
satisfies $\S$ then its predecessor $P(\J)$ satisfies $P(\S)$.
Then Lemma~\ref{recurre1} establishes the inequalities
\begin{equation}\label{recurrenceN}
N(\S)\leq N(P(\S)) + c \times \sum_{i\not\in U}
\sum_{j=1}^d
N(P(\S)\times J_{i,j}^t),
\end{equation}
satisfied for each $(U,t+1)$--specification $\S$, with
$U\subset\{1,2,\ldots,d\}$. We do not solve this chain of
inequalities, but instead we solve the chain corresponding to the
maximal degeneracies. Indeed, for each $U \subset
\{1,2,\ldots,d\}$ and $ t\in \nn $, let
\[
n (U,t) := \max \{ N(\S):\ \S \text{ is a }(U,t)
\text{--specification} \}.
\]
Taking the maximum at both sides of the inequality
(\ref{recurrenceN}), we obtain the chain of recursions
\begin{equation}\label{recurrencen}
n(U,t+1)\leq n(U,t) + c\times \sum_{i\not\in U} (\#\PP_i-1)\times
n(U\cup \{i\},t)
\end{equation}
taking into account that $\#\{1\leq j\leq d: T_{i,j} >
0\}=\#\PP_i-1$. We solved this chain of recursions as follows.

%%%%%%%

\ms Let us first take $U:=\{1,2,\ldots,d\}
\setminus \{i\}$, for $1\leq i\leq d$ arbitrary. Taking into account
that $n(\{1,2,\ldots,d\},t)=1$ for each $t\in\nn$, from (\ref{recurrencen})
we obtain the linear recursion
\[
n(U,t+1)\leq n(U,t) + c\times (\#\PP_i-1).
\]
where $c:=\#\PP-1$.
This recursion has solution
$n(U,t)\leq n(U,1)+ c\times(\#\PP_i-1)\times (t-1)$.
Since $n(U,1)\leq \#\PP_i\leq c\times (\#\PP_i-1)$, we have
\[
n(U,t)\leq c\times(\#\PP_i-1)\times t, \ \ \forall \ t\in\nn.
\]

\ms Suppose that for each
$U\subset \{1,2,\ldots,d\}$ with  
$\#\{1\leq i\leq d:\ i\not\in U\}=k$. 
Then,
\begin{equation}\label{recurrencen2}
n(U,t)\leq c^k\times t^k \times\prod_{i\not\in U}(\#\PP_i-1)\
\ \forall \ t\in\nn.
\end{equation}
Let $U'\subset \{1,2,\ldots,d\}$, with
$\#\{1\leq i\leq d:\ i\not \in U'\}=k+1$. Substituting this 
into %
(\ref{recurrencen}), 
and %
taking into account the induction hypothesis,
we obtain
\begin{eqnarray*}
n(U',t+1)&\leq& n(U',t) + c\times \sum_{i\not\in U'} (\#\PP_i-1)
                      \times n(U'\cup\{i\},t)\\
          &\leq & n(U',t) + (k+1)\ c^{k+1}\times t^k\times
                                \prod_{j\not\in U'}(\#\PP_j-1),
\end{eqnarray*}
with solution $n(U',\tau) \leq n(U',1)+ (k+1)\  c^{k+1}\times
\left(\sum_{t=1}^{\tau-1}t^k \right) \times \prod_{j\not\in
U'}(\#\PP_j-1)$, for each $\tau\in \nn$.

\ms Now, taking into account that $\sum_{t=1}^{\tau-1}t^k\leq
\tau^{k+1}/(k+1)$, we get
\[
n(U',\tau)\leq c^{k+1}\times \tau^{k+1}\times\prod_{j\not\in
U'}(\#\PP_j-1),
\]
which proves inequality (\ref{recurrencen2}) %
for each
$U\subset \{1,2,\ldots,d\}$.

\ms Since $N(J)\leq n(\{j\},t)$ for each $t\in \nn$ and every
$J\in \QQ_{j}^t$,  by %
replacing % 
(\ref{recurrencen2}) into
(\ref{recurrenceC3}), we obtain an autonomous recursion satisfied
by the complexity: 
\begin{eqnarray*}
C(t+1) &\leq & C(t)+ c\times d  c^{d-1}
\times t^{d-1}\times
\prod_{i=1}^d(\#\PP_i-1).
\end{eqnarray*}
We solve this recursion taking into account that $C(1)=\#\PP=c+1$
and that $\prod_{i=1}^d(\#\PP_i-1)\leq c$, and % 
finally obtain the
upper bound
\begin{eqnarray*}
C(\tau)&\leq &1+
        c\left(1+d c^{d}\left(\sum_{t=1}^{\tau-1}t^{d-1}\right)\right)  \\
       &\leq &1+c\left(1+ c^d\ \tau^d\right),
\end{eqnarray*}
for each $\tau\in\nn$.
\endpf

\bs
\subsection{Generalizations} Theorem~\ref{polynomial} admits two quite direct
generalizations as follows.

\ms \paragraph{\bf Rectangular piecewise 
contractions}~\label{general-contractions}
Lemmas from~\ref{injective} to~\ref{recurre1} can be extended 
to piecewise contractions of the following type:
give the base partition $\PP:=\prod_{i=1}^d\PP_i$, let  $F:[0,1]^d\to[0,1]^d$ be
such that for each $1\leq i\leq d$ and $\I\in \PP$, there exists a contractive map
$f_{i,\sI}:[0,1]\to[0,1]$ such that
\[
F(\x)_i=f_{i,\sI}(\x_i), \ \forall \ \x\in\I. 
\]
We say that $F$ so defined is coordinatewise injective when for each 
$1\leq i\leq d$, if $f_{i,\sI}\neq f_{i,\sJ}$ then 
$f_{i,\sI}([0,1])\cap f_{i,\sJ}([0,1])=\emptyset$.
Theorem~\ref{polynomial} holds for systems $([0,1]^d,F)$, where $F$ is
a piecewise contractions of the type just described, satisfying the coordinatewise
injectivity. 

\ms We say that a piecewise contraction $F:[0,1]^d\to[0,1]^d$ of the type 
describe above is a {\it rectangular piecewise contraction} compatible with $\PP$. 
We use the adjective rectangular since such a transformation, when restricted to an 
atom of $\P$, maps the products of intervals into products of intervals.

\ms 
The bound $C(t)\leq 1+c(1+c^dt^d)$ holds for any $d$--dimensional regulatory 
network, as far as the piecewise injectivity is ensured. Hence, 
lemma~\ref{prop-injectivity} implies that given a interaction matrix $M$, all the
regulatory networks with interaction matrix $M$ and sufficiently small contraction rate, 
have complexity bounded by the same polynomial $t\mapsto 1+c(1+c^dt^d)$. 

\ms 
\paragraph{\bf Compatible sequences}~\label{compatible-sequences}
Instead of a single rectangular 
piecewise contraction $F:[0,1]^d\to [0,1]^d$ compatible with a given base
partition $\PP$, we may consider a sequence 
$\SS:=\{F_t:[0,1]^d\to[0,1]^d\}_{t=1}^\infty$ of rectangular piecewise contractions 
compatible with $\PP$. Then, we may consider the non--autonomous
system generated by $\SS$ where orbits are given by 
\[
\OO_{\SS}(\x):=\{\x^{t}:=F_t(\x^{t-1}), \ \x^0=\x\}.
\]
The {\it dynamical complexity of the sequence $\SS$} is the integer function
$C_{\SS}:\nn\to\nn$ such that $C_{\SS}(t):=\#\PP_{\SS}^{t}$, where 
$\PP_{\SS}^1:=\PP$ and 
$\PP_{\SS}^{t+1}:=\{F_{t+1}^{-1}(\I)\cap \J:\ \I\in \PP_{\SS}^t \text{ and } \J\in\PP\}$,
for $t\geq 1$. 
If each one of the maps in $\SS$ satisfies 
the coordinatewise injectivity, then lemmas~\ref{injective} to~\ref{recurre1} hold.
In this case we use
$\QQ_{\SS}^{t}:=F_{t}\circ\cdots\circ F_2\circ F_1(\PP^t_{\SS})$ in place of
$\QQ^t$.  

\ms {\it Mutatis mutandi} one can prove the analogous of Theorem~\ref{polynomial}.

\begin{tmma}~\label{t-sequence}
Let $\SS:=\{F_t:[0,1]^d\to[0,1]^d\}_{t=1}^\infty$ be sequence 
of rectangular piecewise contractions compatible with a partition
$\PP:=\prod_{i=1}^d\PP_i$. If all the contractions in $\SS$ satisfy the
coordinatewise injectivity, then $C_{\SS}(t)\leq 1+c(1+c^dt^d)$
with $c:=\#\PP-1.$
\end{tmma}

\ms In what follows a sequence of rectangular piecewise contractions compatible
with a given partition $\PP$, all of them satisfying the coordinatewise injectivity, 
will be refereed as a {\it $\PP$--compatible sequence}.

\ms
\begin{remark}~\label{various-injective}
The complexity of a compatible sequence of piecewise contractions with base partition
$\PP$ is bounded by the polynomial $t\mapsto 1+c(1+c^dt^d)$. Note that the constant
$c$ depends only on the cardinality of the base partition, indeed, $c:=\#P-1$. 
\end{remark}

\bs
\section{The underlying network}~\label{underlying-network}

\

\ms {\it A network} (directed graph) consist on a finite set $V$ of
{\it vertices}, and a set of ordered couples 
$\AA\subset V\times V$, the {\it arrows}. We consider the interacting
units to form a network with vertex set  $V:=\{1,2,\ldots,d\}$ and arrow set
$A:=\{(i,j)\in V\times V: \ K_{i,j} > 0\}$. We call this {\it the
underlying network}. As we mentioned above, the main motivation of
this study is to understand the relationship between the structure of
the underlying network and the asymptotic dynamics of the system.

\ms Let $T,K:\{1,2,\ldots,d \}\times \{1,2,\ldots,d \}\to [0,1]$ and 
$s:\{1,2,\ldots,d \}\times \{1,2,\ldots,d \}\to \{-1,0,1\}$ be compatible 
matrices. Assume that for each $1\leq j \leq d$, $\sum_{i=1}^dK_{i,j} < 1$ and fix 
a contraction rate $a\in (0,1)$. To the sequence of vectors
$\DD:=\left\{D_t\in [0,1]^d\right\}_{t=1}^\infty$ satisfying 
$0\leq D_{t,j}+\sum_{i\in V}K_{i,j} \leq 1$, 
we associate the sequence 
$\DD_{K,T,s,a}:=\left\{F_t:[0,1]^d\to [0,1]^d\right\}_{t=1}^\infty$ of 
piecewise affine contractions such that 
\[
F_t(\x)_j:=a\x_j+(1-a)\left(\sum_{i\in V}K_{i,j}H(s_{i,j}(T_{i,j}-\x_i))+D_{t,j}\right),
\]
for each $t\in \nn$ and $1\leq i\leq d$.
Notice that threshold matrix $T$ completely determines the partition 
$\PP:=\prod_{i\in V}\PP_i$ whose atoms are the domains of continuity of
the maps $F_t\in\DD_{T,K,s,a}$. Taking $a\in (0,1)$ 
sufficiently small so that each transformation $F_t\in\DD_{T,K,s,a}$ is coordinatewise 
injective, we obtain a $\PP$--compatible sequence. In this way we can define
the complexity
\begin{equation}\label{complexity-bundle}
C_{T,K,s,a}(t):=\max\left\{C_{\DD_{T,K,s,a}}(t):\ 
\DD \text{ such that } 0\leq D_{t,j}+\sum_{i\in V}K_{i,j} \leq 1 \forall t\in\nn\ j\in V\right\}
\end{equation}
for the parameters $K,T,s,a$.

\subsection{Skew product} 
The first result in this direction concerns networks where the dynamics on 
a subnetwork is forced by the dynamics on the complementary subnetwork. 
Let $(V,A)$ be a network such that the vertex set $V$ may be decomposed
into two disjoint sets $V=V_{\rm base}\cup V_{\rm bundle}$ such that
$A\cap (V_{\rm bundle}\times V_{\rm base})=\emptyset$. We say 
in this case that $(A,V)$ admits {\it base---bundle decomposition}. 

\ms
When the underlying network $(V,A)$ of the system $([0,1]^{\#V}, F)$ 
admits a base--bundle decomposition $V=V_b\cup V_d$, 
the associate interaction matrix $K\in\MM_{\#V}([0,1])$ has the triangular 
block form
\[
K=\left(\begin{array}{cc} K_{b,b} & K_{b,d}\\{\bf 0}  & K_{d,d} \end{array}\right),         
\]
where $K_{b,b}:=K|_{V_b\times V_b}$, $K_{b,d}:=K|_{V_b\times V_d}$ and
$K_{d,d}:=K|_{V_d\times V_d}$.  In this case the projection 
$F(\x)_b:=\{\Pi_jF(\x):\ j\in V_b\}$ depends only on 
$\x_b:=\{\Pi_j {\bf x}:\ i\in V_b\}$, thus
the projection of $([0,1]^{\#V},F)$ to coordinates in $V_{\rm base}$ 
is well defined. 

\ms If the underlying network $(V,A)$ of the system $([0,1]^{\#V}, F)$ 
admits a base--bundle decomposition $V=V_b\cup V_d$, then 
the dynamics has the structure of a {\it skew product}. 
Indeed, let 
\[
F_b:[0,1]^{\#V_b}\to[0,1]^{\#V_b} \text{ such that }
F_b(\y):=F(\y\oplus \z)
\]  
with $\z\in [0,1]^{\#V_d}$ arbitrarily chosen, and for each 
$\y\in [0,1]^{\#V_b}$ let 
\[
F_{\sy}:[0,1]^{\#V_b}\to [0,1]^{\#V_b} \text{ such that }
F_{\sy}(\z):=F(\y\oplus\z).\] 
Then, for $\x\equiv\x_b\oplus\x_d\in[0,1]^d$, we have
$F(\x)=F_b(\x_b)\oplus F_{\sx_b}(\x_d)$.

\bs

\begin{tmma}~\label{skew} 
For a regulatory network with parameters $T,K,s,a$, 
suppose that the associated underlying network $(V,A)$ 
admits a base--bundle decomposition $V=V_b\cup V_d$.  
If $([0,1]^{\#V}, F)$ satisfies coordinatewise injectivity, then
\[
C(t)\leq C_b(t)\times C_{(K_d,T_d,s_d,a)}(t), 
\]
where $t\mapsto C_b(t)$ is the dynamical complexity of the regulatory network
$([0,1]^{\#V_b}, F_b)$, and $K_d,T_d$ and $s_d$ denote respectively the restrictions 
of $K,T$ and $s$ to $V_d\times V_d$. 
\end{tmma}

\ms 
\begin{proof}
Let $\PP_b:=\prod_{i\in V_b}\PP_i$ and $\PP_d:=\prod_{i\in V_d}\PP_i$ be 
respectively the projection of the base partition to the coordinates in $V_b$ 
and $V_d$. Because of the skew product structure of the map,
for each $t\in \nn$ and $\I\in \PP^t$ we have intervals $\I_b\in\PP_b^t$ and
$\I_d\in [0,1]^{\#V_d}$ such that $\I=\I_b\times I_d$. Since for each
$\x_b\in [0,1]^{\#V_b}$, the corresponding
bundle transformation $F_{\sx_b}:[0,1]^{\#V}\to [0,1]^{\#V}$ 
depends only on the atom $\J_b\in \PP_b$ containing $\x_b$, 
then for each $\I_b\in \PP_b^t$ there exists a finite sequence of piecewise
transformations $F_1,F_2,\ldots,F_t:[0,1]^{\#V_d}\to [0,1]^{\#V_d}$ such
\[
F^t(\I_b\times \I_d)=F_t\circ F_{t-1}\circ\cdots\circ F_1(\I_d),
\]
for each multidimensional interval $\I_d\subset [0,1]^{\#V_d}$. Now, since $F$ 
is coordinatewise injective, then the finite sequence $F_1,F_2,\ldots, F_t$ can be seen
as the prefix of a $\PP_d$--compatible sequence in the collection 
$\DD^{K_b,T_b,s_d,a}$ for an appropriated sequence of vectors $\DD$.
From this it follows that 
\[
N(\I_b):=\#\{\I\in \PP^t: \I=\I_b\times \I_d, \text{ for some } \I_d\in [0,1]^{\#V_d}\}
\leq C_{K_d,T_d,s_d,a}(t),
\]
for each $\I_b\in \PP^t_b$. From this it follows that 
\[
C(t)\equiv 
\#\PP^t\leq \#\PP_b^t\times C_{K_d,T_d, s_d,a}(t)=C_b(t)\times C_{K_d,T_d, s_d,a}(t),
\]
and the theorem is proved.
\end{proof}

\subsection{Redundant and essential vertices}
Theorem~\ref{polynomial} establishes that the complexity of
a discrete time regulatory network, when the contraction rate is sufficiently
strong, grows at most polynomially. The degree of this polynomial bound
is never larger than the cardinality of the vertex set. In this section we
will see how this degree can be related to the structure of the underlying
network, via the notions of redundant and essential sets of vertices.

\ms For $i, j\in V=\{1,2,\ldots,d\}$ fixed, for each $t\in \nn$,
and for every $J\in \QQ^{t}_{i}$, let
\begin{equation}\label{constrain}
C_j(J)=\Pi_jC(J):=\{\Pi_j \J:\ \J\in \QQ^t\text{ and }\Pi_i\J = J\}.
\end{equation}
In other words, if the $i$--th projection of $\J\in\QQ^t$ is
given, i.~e., if $\J$ satisfies a given $(\{i\},t)$--specification,
and it is equal to $J\in\QQ^t_i$, then its $j$--th projection belongs
to the set $C_j(J)\subset \QQ^t_j$. We say that {\it $i$ drives $j$}
if the cardinality of $C_j(J)$ is uniformly bounded as $J$ ranges over
$\cup_{t=1}^\infty\QQ_i^t$, i. e., if there exists
a constant $m_{i,j}\in \nn$ such that
$\#C_j(J)\leq m_{i,j}$ for all $J\in \cup_{t=1}^\infty\QQ_i^t$.

\ms A set of vertices $U\subset V$ is {\it redundant} if for each
$i\in U$ there exists $j\in V\setminus U$ such that $i$ drives $j$
in the sense stated above. Below we will see why the driving vertex
is the redundant one. 

\ms On the other hand, a set $W\subset V$ is {\it essential} if
there exists a fixed multiplicity $M\in \nn$ such that for each
$t\in \nn$ and $\S\in \prod_{i\in W}\QQ_i^t$,
$N(\S)\leq M$, i.~e. $W\subset V$ is essential if the degeneracy 
of the $(W,t)$--specifications is uniformly bounded as $t$ varies. 
Of course $W=V$ is essential, but the interesting 
case arises when there is a essential proper subset of the vertex 
set, whose complement is redundant.

\ms
\begin{tmma}[Bounding the degree]\label{bound-degree}
Let $([0,1]^d, F)$ be a discrete time regulatory network
satisfying coordinatewise injectivity, and let $(V,\AA)$ be the
underlying network. Suppose that $U\subset V$ is redundant,
whereas $W:=V\setminus U$ is essential. Then, there are constants
constant $c_1,c_2 > 0$ such that for each $t\in \nn$
$C(t)\leq 1+c_1\left(1+c_2t^q\right)$, with $q:=\#W$.
\end{tmma}

\ms The proof of this theorem follows the same scheme, with the 
suitable modifications, as the proof of Theorem~\ref{polynomial}.

%%%%%%%%%%%%%%%%%%%%%%%%%%%%%%
\begin{proof} Since we are assuming coordinatewise injectivity, 
Lemma \ref{injective} ensures that $F$ is injective, so that $C(t)=\#\QQ^t$.
As in the proof of Theorem~\ref{polynomial}, we use 
the recursion in~(\ref{recurrenceC}), which can be written as
\[
C(t+1)\leq C(t)+c_1\times \#\{\J\in\QQ^{t}:\ \Pi_i\J=J^{t}_{i,j}
\text{ for some } (i,j)\in\AA\},
\]
with $J^t_{i,j}$ as defined by Corollary~\ref{determine}.

\ms From here on we use the notion of redundant and essential
vertices. Fix $U\subset V$ redundant, such that 
$W:=V\setminus U$ is essential. Now, since each essential vertex 
can be driven by several redundant ones, let us fix a partition 
$\CC:=\{R_j\subset U:\ j\in W\}$ of $U$, 
where to each atom of the partition therecorresponds an essential vertex, 
such that $i$ drives $j$ for each $i\in R_j$. For the rest of the
essential vertices let $R_j=\emptyset$. 
Such partition is not unique since different essential
vertices may driven by the same redundant one. For $\CC$ fixed,
for each $t\in \nn$ and $j\in W$ define the collection of intervals
\begin{equation}\label{redundant}
\RR_j^t:=\bigcup_{i\in R_j}\bigcup_{(i, k)\in \AA} C_i(J_{i, k}^t),
\end{equation}
with $C_i(J_{i,k}^t)$ as defined by (\ref{constrain}).

%%%%%%%%

\ms The vertices in $R_j$ are redundant and $\#\{k\in V:\
(i, k)\in\AA \}\leq \PP_i-1$. Therefore
\begin{equation}\label{redundancy}
\#\RR_j^t\leq m_j:=(\#\PP_j-1)+\sum_{i\in R_j} m_{i, j}\times (\#\PP_i-1)
\end{equation}
for each $t\in \nn$. Here $m_{i,j}$ is the multiplicity
bounding the cardinality of $C_i(J_{i,j}^t)$. For a fixed partition $\CC$ and
taking into account inequality (\ref{recurrenceC}), we obtain 
\begin{equation}\label{recurrenceC4}
C(t+1)\leq C(t)+c_1\times 
\sum_{j\in W}\sum_{J\in\RR_{j}^t} N(J),
\end{equation}
Let us recall that $N(J):=\#C(J)=\#\{\J\in\QQ^t:\ \Pi_i\J=J\}$. 

\ms As in the proof of Theorem~\ref{polynomial}, we can find 
polynomial bounds for the growth of the degeneracies $N(J)$. 
Taking into account the fact that redundant vertices drive 
essential ones, we may rewrite inequality~(\ref{recurrenceN}) as
\[
N(\S)\leq N(P(\S)) + c_1 \times \sum_{j\in W\setminus
W'}\sum_{J\in \RR_j^t} N(P(\S)\times J),
\]
satisfied for each $(W',t+1)$--specification $\S$, with $W'\subset
W$. Once again we solve the chain of inequalities satisfied by 
the maximal degeneracies,
\begin{equation}\label{recurrencenR}
n(W',t+1)\leq n(V',t) + c_1\times \sum_{j\in W\setminus W'}m_j\
n(W'\cup \{j\},t).
\end{equation}
Using exactly the same scheme as in the proof of Theorem~\ref{polynomial} 
we obtain
\begin{equation}\label{recurrencen2R}
n(W',t)\leq M\ c_1^k\ t^k \ \prod_{j\in W\setminus W'} n_j, \ \
\forall \ t\in\nn,
\end{equation} 
with $n_j:=\max(m_j,\#\PP_j)$. 
%%%%%%%%%%%%%

\ms Since $N(J)\leq n(\{j\},t)$ for each $t\in \nn$ and every
$J\in \QQ_{j}^t$,  substituting inequality (\ref{recurrencen2R}) into
(\ref{recurrenceC4}), we obtain 
\begin{eqnarray*}
C(t+1) &\leq & C(t)+ c_1\ \sum_{j\in W} \sum_{J\in \RR_j^t}M\
c_1^{q-1}\ t^{q-1}\prod_{i\in W\setminus\{j\}}n_i\\ &\leq & C(t)+
M\ c_1^{q}\ t^{q-1} \ \sum_{j\in W}m_j \prod_{i\in W\setminus
\{j\}}n_j,
\end{eqnarray*}
with $q=\#W$. We solve this recursion by taking into account that
$n_j\geq m_j$ and $C(1)=c_1+1$, to finally obtain the upper bound
\begin{eqnarray*}
C(\tau)&\leq &1+c_1\left(1+q\ M\ c_1^{q}\prod_{j\in W}n_j
\left(\sum_{t=1}^{\tau-1}t^{q-1}\right)\right)  \\
       &\leq &1+c_1\left(1+ M\ c_1^q\ \tau^q \prod_{j\in W}n_j \right),
\end{eqnarray*}
for each $\tau\in\nn$. The theorem follows by making $c_2=M\ c_1^q\
\prod_{j\in W}n_j$.
\end{proof}

%%%%%%%%%%%%%%%%%%%%%%%%%%%%%%%%%%%

\bs 

\subsection{Head--independent sets} Now we will 
exhibit concrete situations where the underlying 
network $(V,\AA)$ admits a non--trivial redundant vertex set with 
essential complement. In this way we obtain concrete examples 
of discrete time regulatory networks whose dynamical complexity 
is bounded by a polynomial of degree strictly smaller than $\#V$.

\ms Let us first introduce a characterization of a
class of vertices which are essential whenever coordinatewise
injectivity holds.

\ms Fix a network $(V,\AA)$. A set of vertices $U\subset
V:=\{1,2,\ldots,d\}$ is {\it head--independent} if
\begin{itemize}
\item[1)] $\{j\in V:(i,j)\in \AA\}\cap U=\emptyset$
for all $i\in U$, and 
\item[2)]
$\{j\in V:\ (i,j)\in \AA\}\cap\{j\in V:\ (i',j)\in
\AA\}=\emptyset$ for all $i,i'\in U$ with $i\neq i'$.
\end{itemize}
In other words, $U\subset V$ is head--independent if 1) each arrow
with tail in $U$ has head outside $U$, and 2) different arrows
with tail in $U$ have necessarily different heads. We will show
that for discrete time regulatory network satisfying coordinatewise
injectivity, the complement of a head--independent set is
essential.

\ms
%%%%%%%%%%%%%%%%%%%%%%%%%%%%%%%%%%%%%%

\begin{proposition}\label{complement-head}
Suppose that the discrete time regulatory network $([0,1]^d,F)$ 
satisfies coordinatewise injectivity. Let $(V,\AA)$ be the underlying
network, let $U\subset V$ be head--independent, and let $W:=V\setminus U$.
Then for each $\I\subset \PP$, each $t\in \nn$, and every
$(W,t)$--specification $\S$ we have
\[ \#\{\J\in C(\S):\ \J\subset \I\}\leq 1.\]
\end{proposition}

\begin{proof}
The statement holds for $t=1$ since in this case if $\J\in
\QQ^1\equiv\PP$ is such $\J\subset \I$, then necessarily $\J=\I$.
Now, given a $(W,1)$--specification $\S$ there are two
possibilities: 1) either $\Pi_i\S=\Pi_i\I$ for each $i\in W$, or
2) there exists $i\in W$ such that $\Pi_i\S\cap\Pi_i\I=\emptyset$.
In the first case $\#\{\J\in C(\S):\ \J\subset \I \}=1$, while in
the second $\#\{\J\in C(\S):\ \J\subset \I\}=0$.

\ms Let us suppose that $\#\{\J\in C(\R):\ \J\subset \I\}\leq 1$
for $\tau\in\nn$, for each $(W,\tau)$--specification $\R$, and for
every atom $\I\in \PP$. Fix a $(W,\tau+1)$--specification $\S$.
According to Lemma \ref{predecessor}, if $\J\in \QQ^{\tau+1}$
satisfies $\S$ then  $P(\J)$ satisfies $P(\S)$. Therefore
\begin{eqnarray*}
\{\J\in C(\S):\ \J\subset \I\} &\subset & \{F(\J')\cap \I:\ \J'\in
             C(P(\S)) \}\\
                                 &=&\cup_{\sI'\in \PP}
\{F(\J')\cap\I:\ \J'\in C(P(\S))\text{ and } \J'\subset \I'\}.
\end{eqnarray*}
By hypothesis $\#\{F(\J')\cap \I:\ \J'\in C(P(\S))\text{ and }
\J'\subset \I'\}\leq 1$ for each $\I'\in\PP$. Now, for $\I'\in
\PP$ to be such that $\{F(\J')\cap \I:\ \J'\in C(P(\S))\text{ and
}\J'\subset\I'\}\neq\emptyset$, it is needed that
$\Pi_j\I'\supset\Pi_jP(\S)$ for each $j\in W$, which already
determines $\Pi_j\I'$ for $j\in W$. At the same time, for
$\J'\subset\I'$ and $j\in W$ we have $\Pi_j(F(\J')\cap \I) =
\S_j$. Therefore the system of constraints
\begin{equation}\label{head-independent}
\left(a\Pi_j\J'+(1-a)\sum_{i=1}^dK_{i,j}H(s_{i,j}(\x_i -
T_{i,j}))\right)\cap\Pi_j \I= \Pi_j\S,\ \forall\ j\in W,
\end{equation}
must hold for any $\x\in \J'$. We consider this system of
constraints as a system of equations for the unknown binary values
$\{b_{i,j}:=H(s_{i,j}(\x_i-T_{i,j})):\ i\in U, \ (i,j)\in \AA\}$.
Since $U$ is head--independent, each one of those equations
contains at most one of these unknowns, and each unknown appears
in at least one of those equations. Coordinatewise injectivity
implies that for $j\in W$ fixed, the two possible values for
$y_{i,j}$ appearing in (\ref{head-independent}) determine two
disjoint intervals. Hence, if $b_{i,j}=b$ is a solution for
(\ref{head-independent}), then $b_{i,j}=1-b$ is not. Thus, the system
of constraints in (\ref{head-independent}) cannot be satisfied by
two different choices $\{b_{i,j}\in \{0,1\}:\ i\in U, \ (i,j)\in
\AA\}$. Finally, since each projection $\Pi_i\I'$ is completely
determined by the values $\{H(s_{i,j}(\x_i-T_{i,j})):\ (i,j)\in
\AA, \x_i\in \Pi_i\J'\}$, then the projections $\Pi_i\J', \ i\in
U$ are determined as well. In this way we fix the remaining
projection $\Pi_i\I'\supset \Pi_i\J'$, proving that there is at
most one atom $\I'\in\PP$ such that $\#\{\J'\in C(P(\S)):\
\J'\subset \I',\ F\J'\cap \I\neq \emptyset\}=1$.
\end{proof}

\ms
\begin{corollary}[Complement of a head--independent is 
essential]\label{c-head-independent} Under the same hypotheses as in the
previous proposition, the set $W:=V\setminus U$ is essential.
\end{corollary}

\begin{proof}
Indeed, since for each $I\subset \PP$, each $t\in \nn$, and every
$(W,t)$--specification $\S$ we have $\#\{\J\in C(\S):\ \J \subset
\I\} \leq 1$, then
$N(\S)=\sum_{\sI\in P}\#\{\J\in C(\S):\ \J \subset \I\}\leq \#\PP$,
which proves the claimed result.
\end{proof}

\ms
\begin{remark}\label{rm-c-head-independet}
The previous corollary holds if instead of the iterates of a
coordinatewise injective map we consider a $\PP$--compatible sequence of 
$\DD_{K,T,s,a}$, such that $U\subset V$ is head--independent 
with respect to the underlying network $(V,A)$ associated to $K$.
\end{remark}

%%%%%%%%%%%%%%%%%%%%%%%%%%%

\bs 

\subsection{2--loops}~\label{ss-2-loops} In order to have an 
effective reduction of 
degree in the
polynomial bound for the complexity, we need that the underlying
network admits redundant sets with essential complement. In
general, head--independent sets are not redundant, and it seems
more difficult to find sufficient conditions for vertices to be
redundant than sufficient conditions for vertex sets to be 
essential.
Below we will show a case where the structure of a two--vertices
subnetwork is such that one of its vertices drives the other.

\ms We say that the interaction matrix $K\in \MM_{d\times
d}([0,1])$ is {\it non--degenerated} if for each $j\in
\{1,2,\ldots,d\}$, the map
\[\left\{\e\in \{0,1\}^d:\ \e_i=0\text{ for all } i \text{ such that }
K_{i,j}=0\right\} \ni \e\mapsto \sum_{i=1}^d\e_iK_{i,j}\in [0,1]\]
is one--to--one.

\ms Fix a network $(V,\AA)$. {\it A 2--loop in $(V,\AA)$} is a couple
$\{i,j\}\subset V$ such that $(i,j)\in \AA$ and $i$ is the only
tail for an arrow whose head is $j$, and $(j,i)\in \AA$ and $i$ is
the only head for an arrow whose tail is $j$. In this case we say
that {\it $j$ is the isolated end of the 2--loop}.

%%%%%%%%%%%%%%%%%%%%%%%%%

\begin{proposition}\label{loop}
Suppose that $([0,1]^d,F)$ satisfies coordinatewise
injectivity, and suppose that $\{i,j\}\subset V$ is a 2--loop in the
underlying network $(V,\AA)$ whose isolated end is the vertex $j$.
If in addition the interaction matrix $K$ is non--degenerated,
then $i$ drives $j$.
\end{proposition}

\begin{proof}
We will prove that, for each $t\in \nn$, each $J\in \QQ^t_i$, and
every $I\in \PP_j$, there exists at most one interval
$J'\in\QQ^t_j$ such that $J'\subset I$, and $C(J\times J')\neq
\emptyset$. In this way we prove that $i$ drives $j$ with
multiplicity  $m_{i,j}\leq \#\PP_j=2$.

\ms For $t=1$ the statement obviously holds. Suppose it holds for
$t=\tau$. Fix $J\in \QQ^{\tau+1}_i$ and $I\subset \PP_j$. Then,
thanks to the coordinatewise injectivity, there exists a unique
contraction $f\in \FF_i$ and a unique interval $P(J)\in
\QQ_j^\tau$ such that $J=f(P(J))\cap I'$. Here $I'$ is the atom
$\PP_{i}$ which contains $J$.

\ms The affine contraction $f:[0,1]\to [0,1]$ has the form
$f(x)=ax+(1-a)\sum_{k=1}^d\e_kK_{k,i}$, for some $\e\in\{0,1\}^d$.
Since the interaction matrix is non--degenerated, then the values
of $\e_k\in\{0,1\}$, for $K_{k,i} > 0$, are uniquely determined by
$f$. Therefore, given $J\in\QQ_i^{\tau+1}$, if $\J'=P(\J)$ for
some $\J\in C(J)$, not only we have $\J'\in C(P(J))$ as
established in Lemma \ref{predecessor}, but also
$H(s_{k,i}(\x_k-T_{k,i}))=\e_k$, for each $\x\in\J'$.

\ms From here on we use the fact that $\{i,j\}$ is a 2--loop with
isolated end $j$, i.~e., (a) $K_{j,k}=0$ for $k\neq i$, and (b)
$K_{k,j}=0$ for $k\neq i$. Condition (a) implies that the value
$\e_j\in \{0,1\}$ determines the atom $I''\in \PP_j$ containing
$\Pi_j\J'$ for each $\J'\in C(P(J))$. Here we apply the induction
hypothesis. Given $P(J)\in \QQ_i^{\tau}$ and $I''\in \PP_j$, there
is at most one interval in $\QQ_j^{\tau}$, which we denote $J'$,
such that $J'\subset I''$ and $C(P(J)\times J')\neq \emptyset$. By
condition (b), the image of each $\J'\in C(P(J)\times J')$ is the
interval $\Pi_jF(\J')=aJ'+(1-a)H(s_{i,j}(\x_i-T_{i,j}))$, which
depends only on $J$ and $P(J)\ni\x_i$. Thus, the interval
$\Pi_jF(\J')\cap I\in \QQ^{\tau+1}_j$ is uniquely determined from
$J\in \QQ^{\tau+1}_i$ and $I\in\PP_j$.
\end{proof}

\ms
\begin{remark}
Once again, proposition~\ref{loop} holds for compatible sequence of 
maps with common base partition, all of the satisfying piecewise injectivity,
and having the same underlying network.
\end{remark}

\bs

\section{Examples}~\label{s-examples}

\subsection{\bf Circuits} 
The circuit on $d$ vertices has vertex set
$V=\{1,2,\ldots, d\}$ and arrow set $\AA\subset V\times V$
satisfying
\[
(i,j)\in \AA \Leftrightarrow \ j=i+1, \text{ or } i=d \text{ and }
j=1.
\]

\ms Some numerical experiments suggest that in large random 
networks, and for a large class of initial conditions, the state of most of the
units in the network converges to a fixed value, whereas in a small
subnetwork, an oscillatory regime takes place. Very often this 
subnetwork is a circuit of low dimension. The circuit on $d$ vertices 
may be negative or positive according to the value of the product
$\sigma:=s_{d,1}\times\prod_{i=1}^{d-1}s_{i,i+1}$, defined by the
activation matrix. In discrete time regulatory networks 
defined in positive circuits, numerical experiments have shown 
mutistability, i.~e., coexistence of multiple attractive fixed points.
On the contrary, for the systems defined on negative circuits, beside 
some degenerated cases, the dynamics does not admit fixed points. 
A rigorous study of this numerically observed phenomena is in
progress. Concerning the complexity of this kind of systems we can
already state the following.

\ms
\begin{center}
\includegraphics[width=0.3\textwidth]{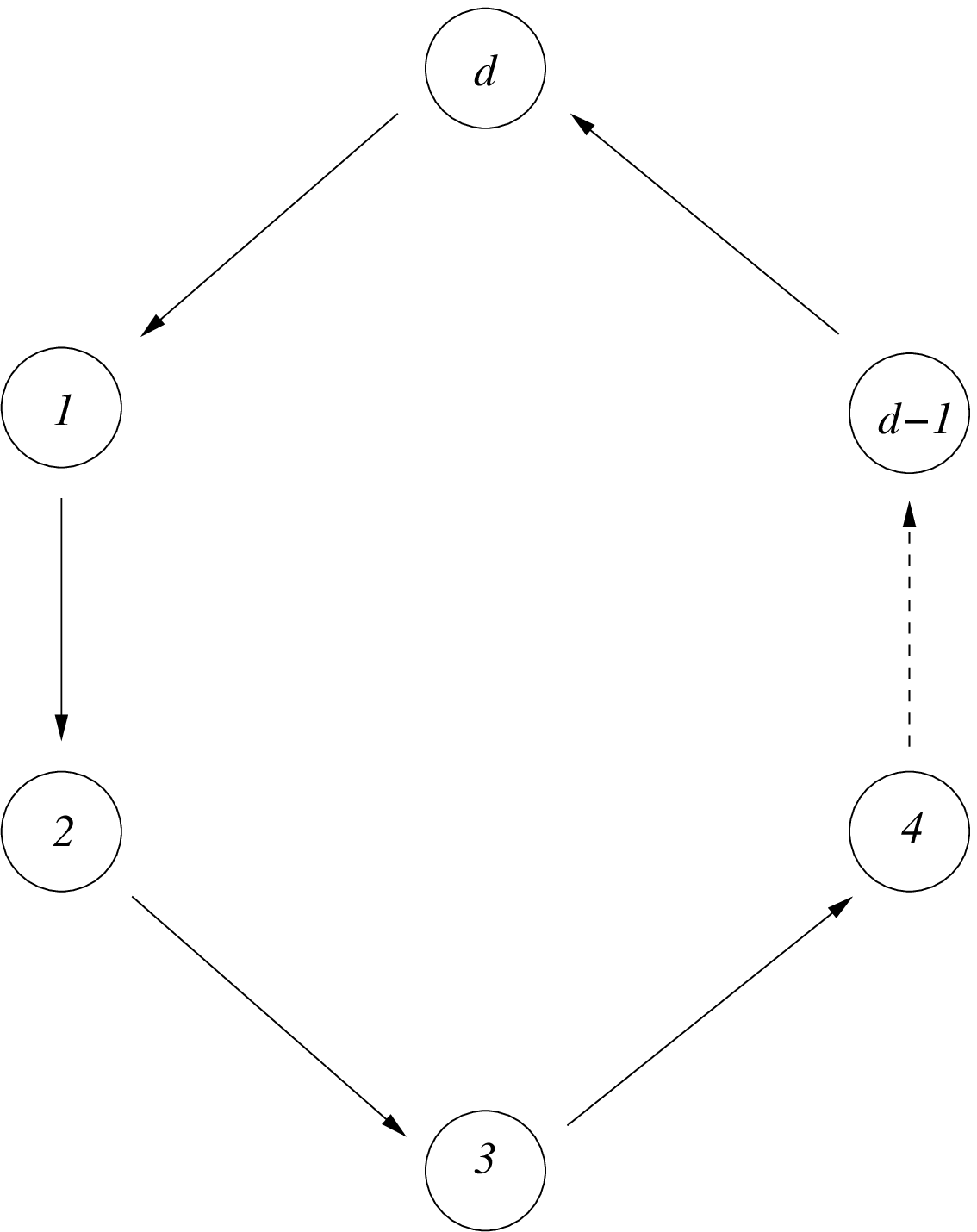}\\
{\bf Figure 1.}
\end{center}

\ms A discrete time regulatory network whose underlying
network is a circuit on $d$ vertices, and whose contraction rate
is smaller than $1/2$, satisfies $C(t)=\OO(t^d)$. This follows from
Theorem~\ref{polynomial}, taking into account that these systems
satisfy coordinatewise injectivity. Indeed, since the
interaction matrix is normalized, i.~e., $\sum_{i=1}^dK_{i,j}=1$,
then $K_{i,j}=1$ for each $(i,j)\in \AA$. Using the same
notation as in Proposition~\ref{prop-injectivity}, we have that 
in this case $\II_i=\{0,1\}$ for each $i\in\{1,2,\ldots,d\}$,
implying that $a_0=1/2$.

\ms The genetic toggle switch described in~\cite{Gardner2000} is
a concrete experimental example of a positive 2--circuit, made of
two mutually inhibiting interacting units. In this case each unit is composed
by a couple promoter--repressor, such that each promoter is inhibited by
the repressor transcribed by the opposite promoter. This system is
traditionally modeled by the system of ordinary differential equations
\begin{eqnarray*}
\frac{dr_1}{dt}=\frac{\alpha_1}{1+r_2^{\beta_2}}-r_1\\
\frac{dr_2}{dt}=\frac{\alpha_2}{1+r_1^{\beta_1}}-r_2,
\end{eqnarray*}
where $r_i$ represents the concentration of repressor $i$,
$\alpha_i$ is the effective rate of synthesis of repressor $i$, and $\beta_i$ is
the cooperativity of repression of promoter $i$. The corresponding
network of piecewise contractions has the form
\begin{eqnarray*}
r_1(t+1)=ar_1(t)+(1-a)H(T_{1,2}-r_2(t))\\
r_2(t+1)=ar_2(t)+(1-a)H(T_{2,1}-r_1(t)),
\end{eqnarray*}
where the rate of synthesis of the repressors is controlled by 
the parameter $a$, and the cooperativity of repression of promoters 
is related to both $a$ and the thresholds $T_{i,j}$.  The dynamical behavior
of the positive circuit in two vertices is completely understood (see~\cite{bastien}),
and it can be either bistable (the same as the continuous model) or 
conjugated to a rotation. In all cases the dynamical complexity asymptotically 
bounded by $t\mapsto t+1$.

\ms For the general circuit on 2 vertices, if the
contraction rate is smaller that $1/2$, each one of the 2 vertices
is at the same time redundant and essential. Indeed, in this case
the coordinatewise injectivity holds, each vertex constitutes a
head--independent set, and at the same time it is the isolated end of 
a 2--loop. Corollary \ref{c-head-independent} and Proposition~\ref{loop} 
apply, implying that each vertex is redundant and essential at the same
time. Hence, a system of two interacting units whose underlying
network is a circuit on 2 vertices, and whose contraction rate is
smaller than 1/2, necessarily has linear complexity. Furthermore,
for the negative circuit, since the image of rectangle in
$\QQ^\tau$ cannot intersect the double discontinuity
$(\x_1=T_{1,2}, \x_2=T_{2,1})$, then it cannot have more than one
successor. Taking into account this constraint in the computations
we performed above, we deduce that in this case $C(t)\leq 2t+2$
for each $t\in \nn$. It has been proved~\cite{bastien} that the
negative circuit admits an arbitrary large number of periodic
orbits, as one moves on the line $\{(T,a):\ T_{1,2}=T_{2,1}=1/2,
a\in [1/2,1)\}$. There the number of admissible periodic orbits grow
as one increases the contraction rate. Because of this, the
complexity increases as $a$ grows, though our numerical
experiments indicate that $C(t)$ can always be bounded by a
quadratic polynomial, as it is illustrated in Figure~2.
We conjecture that the complexity of the
negative circuit on 2 vertices can always be bounded by a linear
function $C(t)\leq K\times (t+1)$, where the constant factor
depends on the choice of parameters $(T,a)$, and it diverges as we
increase $a$.

\ms

\begin{center}
\includegraphics[width=0.5\textwidth]{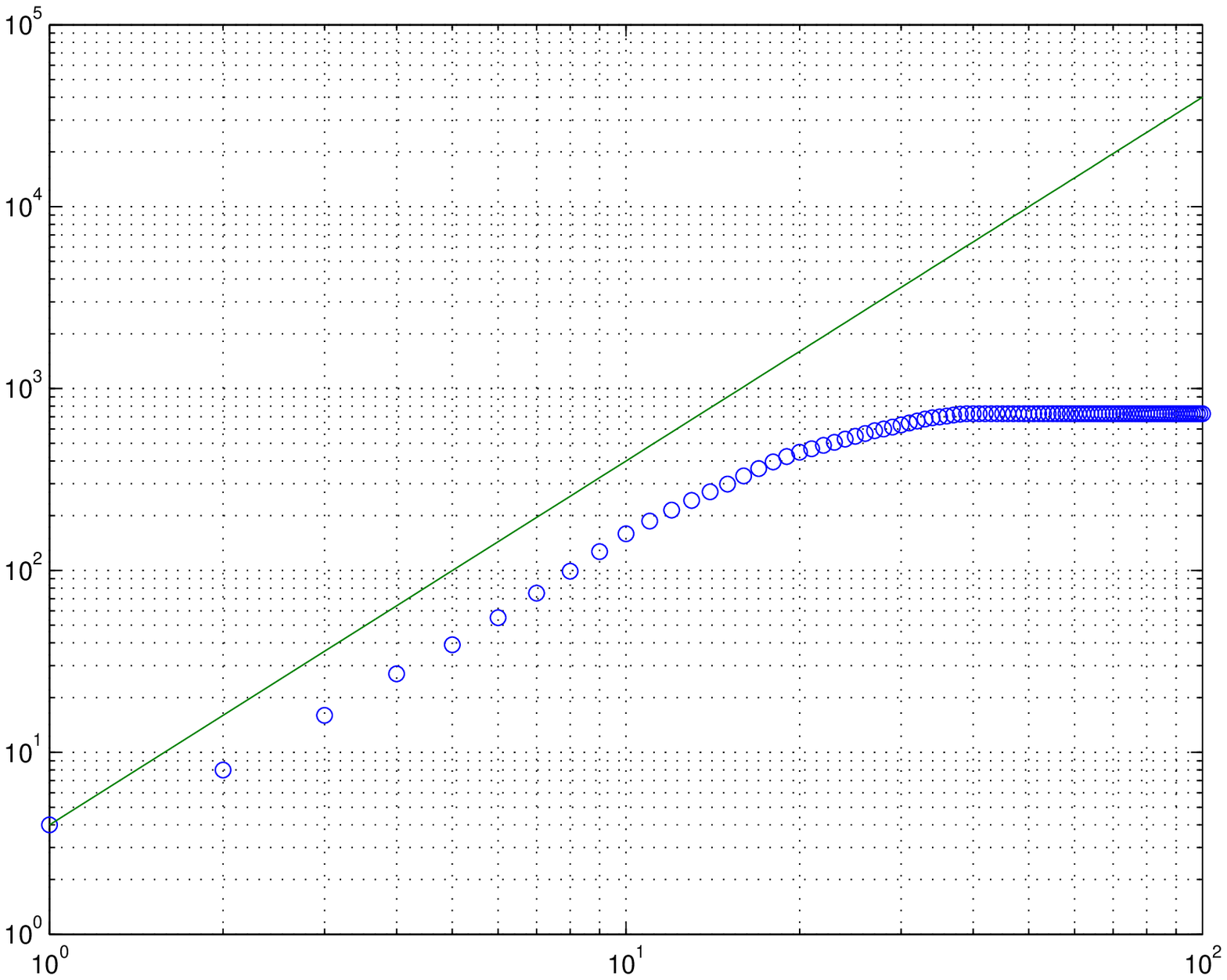}\\
{\bf Figure 2.}\\

\ms \begin{minipage}{0.7\textwidth} Log--log plot of the
complexity of the negative circuit on 2 vertices for
$T_{1,2}=T_{2,1}=0.5$ and $a=0.93$. The straight line is the
log--log plot of the quadratic upper bound $t\to 4 t^2$.
\end{minipage}
\end{center}

\ms For a circuit on $d$ vertices any independent set of vertices
(i.~e., no two of them are adjacent to the same arrow) is
head--independent. In particular
$\{2,4,\ldots, 2\lfloor d/2 \rfloor\}$ is a maximal independent
set, therefore, according to Corollary~\ref{c-head-independent} 
its complement $\{1,3,\ldots, 2\lfloor (d-1)/2 \rfloor+1\}$ is essential.  
If we think of the collection $\QQ^\tau$ as a finite time version of 
the attractor, Proposition~\ref{complement-head} establishes that this coarse 
grained attractor has dimension $k=\lfloor d/2\rfloor$. Indeed, in
this case each $d$--dimensional interval of $\J\in\QQ^\tau$ has at
most $k$ independent one--dimensional projections, whereas $d-k$
projections are determined by the first, and by the atom of $\PP$
containing $\J$. In this way, the coarse grained version of the
attractor lies inside the graph of a $\#\PP$--valued function
defined on $[0,1]^k$, and taking values on $[0,1]^{d-k}$. Note
that this property holds uniformly on $\tau$, suggesting that the
attractor itself may be embedded in the graph of a multivalued
function. This of course has yet to be proved.

\ms Let us remark that in~\cite{Elowitz2000} it is presented an biological 
example of a negative circuit on 3 vertices, the so--called repressiliator. 
They report the existence of oscillatory behavior in time, which would correspond
to periodic solutions of the corresponding model, though the period and amplitude
exhibit significant variability. This may correspond either to stochastic effects, 
or may reflect intrinsic complex behavior.

\subsection{Networks with loops} 
With respect to the effective reduction of the polynomial
upper bound, the example of the circuit on 2 vertices with
contraction rate smaller than $1/2$ may be generalized as follows.
Fix a network $(V,\AA)$ and 2--loops $\{i,j\}$ and $\{i',j'\}$, with
isolated ends $j$ and $j'$ respectively. We say that ``$\{i,j\}$
and $\{i',j'\}$ are disjoint'' if $\{k\in V: \ (i,k)\in \AA\}\cap
\{k\in V:\ (i',k)\in \AA\}=\emptyset$, i.~e., the non--isolated
ends have disjoint heads.

\ms Let $(F,[0,1]^d)$ be a network with non--degenerated interaction
matrix (as defined in the previous section), and satisfying the
coordinatewise injectivity. If the underlying network $(V,\AA)$
has $k:=d$ partwise disjoint 2--loops, then $C(t)=\OO(t^{d-k})$.
Indeed, let $\{\{i_n,j_n\}:\ n=1,2,\ldots, k\}$ be the set of
mutually disjoint 2--loops, and let $j_n$ is the isolated end of
$\{i_n,j_n\}$ for each $n=1,2,\ldots,k$. According to Proposition
\ref{loop}, the set $U:=\{i_n:\ n=1,2,\ldots,k\}$ is redundant.
Indeed, each $j_n$ drives $i_n$ for each $n=1,2,\ldots,k$. On the
other hand, since the 2--loops are partiwise disjoint, the set $U$
head--independent too, and so by Corollary
\ref{c-head-independent} the complement $W:=V\setminus U$ is
essential, therefore $C(t)=\OO(t^{d-k})$.

%%%%%%%%%%%%%%%%%%%%%
\ms
\begin{center}
\includegraphics[width=0.65\textwidth]{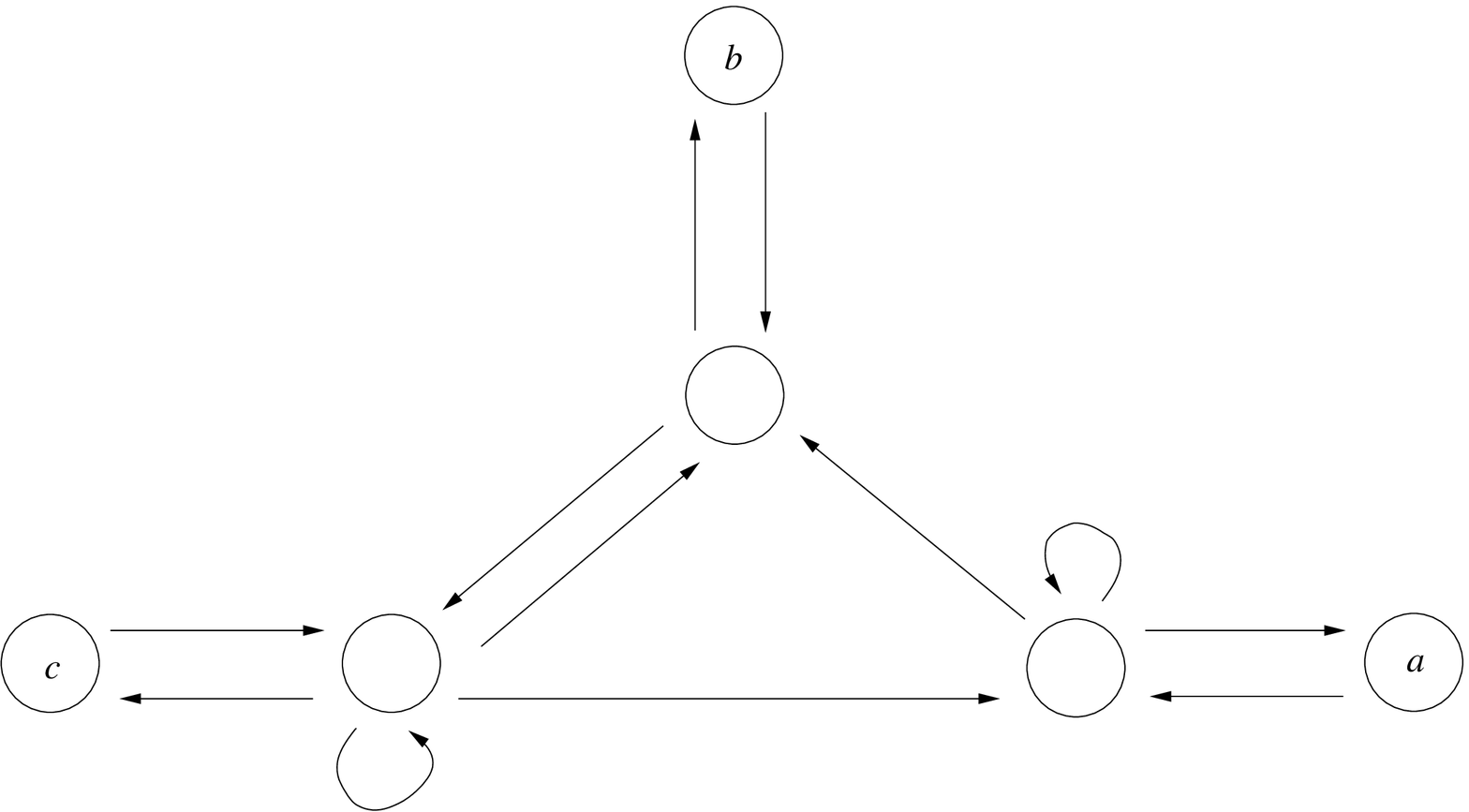}\\
\ms {\bf Figure 3.}
\end{center}

\ms In Figure 3 we represent a network containing 3 disjoint 2--loops
with isolated ends a, b and c. According to
Theorem~\ref{polynomial} and Corollary~\ref{loop}, if the
contraction rate is sufficiently small and the interaction matrix
is non--degenerated, then the complexity of the system is bounded
by a cubic polynomial.

\ms In~\cite{Vilela2004} a discrete time regulatory network of 
the type considered in this paper was used to model the regulation 
of the expression of the tumor suppressor gene {\it p53}, whose product
acts as an inhibitor of uncontrolled cell growth.  
The underlying network of this model is shown in Figure 4. The node {\it m} represents
the expression level of the {\it mdm2} gene, whose product inhibits the expression of {\it p53}.
Nodes {\it b} and {\it c} represent blood--vessel formation and cell proliferation respectively.
Though {\it b} and {\it c} are not expression levels of a particular gene, they represent
relevant quantities directly implicated in cancer decease.

\begin{center}
\includegraphics[width=0.35\textwidth]{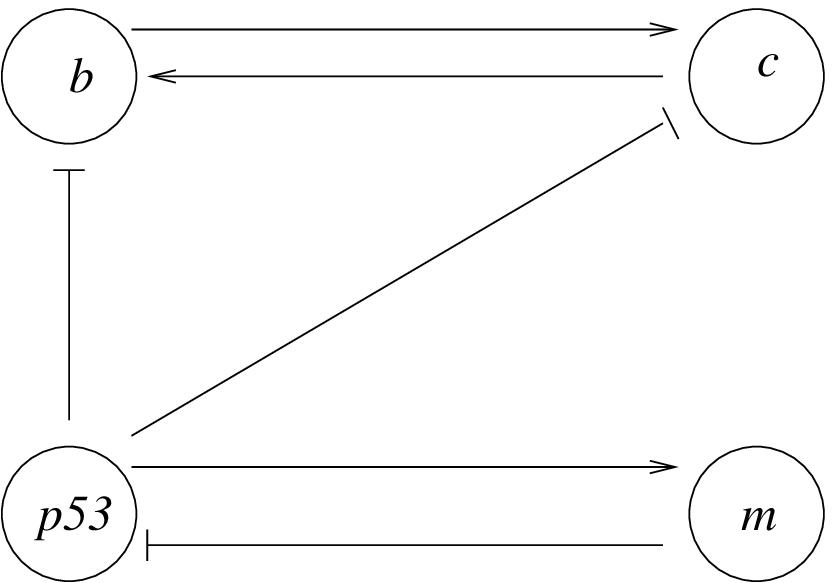}\\
\ms {\bf Figure 4.}\\
\ms \begin{minipage}{0.7\textwidth} Regulatory network for the 
{\it p53} gene expression. The arrow indicates positive interaction, while 
$\dashv$ indicates an inhibitory interaction. 
\end{minipage}
\end{center}

\ms In this network the couple $p53$--$m$ form a 2--loop with 
isolated node $m$, and the vertices $m$ and $c$ are  head--independent. 
If we assume coordinatewise injectivity, according to
proposition~\ref{complement-head}, nodes $p53$ and $c$ 
capture the whole complexity of the network. On the other hand,  
since node $m$ is redundant, the complexity of the network is bounded
by a polynomial of degree 3. Note that in this case the network admits the
base--bundle decomposition $\{p53,m\}\cup\{b,c\}$, and both base ($V_b:=\{p53,m\}$) 
and bundle ($V_d:=\{b,c\}$) define circuit on two vertices. Therefore, if
coordinatewise injectivity holds,
Theorem~\ref{skew} implies that the complexity of this four--vertices network
is bounded by a quadratic polynomial.

\bs
\section{Final Remarks}
\

\ms In this work we formulate and show the first results towards the solution of the
following theoretical problem: determine to what extent the topological structure 
of the underlying network determines the behavior of a discrete time regulatory
network. By solving this problem we intend to acquire a better understanding of the 
behavior observed in real life regulatory networks. The next step in our program consist
in extending Theorem~\ref{polynomial} to the case of 
arbitrary contraction rates, when coordinatewise injectivity does not hold. Some preliminary 
numerical experiments support lead us to conjecture that the complexity is always bounded
by a polynomial; furthermore, this experiments suggest that the asymptotic dynamics 
consist of a finite number of periodic or quasiperiodic orbits. If this is true, the complexity 
should be always bounded by a polynomial of degree one. 

\ms Though in general it is not possible to deduce relevant dynamical properties 
only from the complexity, we expect that in our case the growth of the complexity 
could give information about topological and recurrence properties of the attractor.
For instances, the growth rate of the complexity could be related to the number of 
transitive components of the attractor. Anyway, in order to have a complete description 
of the qualitative behavior of our models we need to consider other characteristics 
besides the dynamical complexity. 

\ms Theorems~\ref{bound-degree} and~\ref{skew} we have specific instances of constraints 
imposed on the behaviors of the system by the structure of the underlying network. 
In this direction,
our next goal is to study how the behavior of a network is related to the possible dynamical
behaviors of its subnetworks. Results in this direction would allow us to deal with  
questions such as the role of self--regulation in the stabilization of the circuits,
and in general, to progress in the understanding of the engineering behind the design of
real life functional networks. 

\bs

\end{document}